\documentclass[final,3p,times,authoryear]{elsarticle}
\usepackage{amsmath, amssymb, amsfonts}
\usepackage{natbib}
\biboptions{authoryear,round}
\usepackage{lmodern}
\usepackage{float}
\usepackage{graphicx}
\usepackage{booktabs}
\usepackage{xcolor}
\usepackage{caption}
\usepackage{subcaption}
\usepackage[utf8]{inputenc}
\usepackage{kotex}
\usepackage{tikz}
\usetikzlibrary{shapes.geometric, arrows.meta, positioning, calc, chains, decorations.pathreplacing, matrix, fit, backgrounds}

\tikzset{
    base/.style={draw, thick, align=center, minimum height=0.8cm, inner sep=4pt, font=\small},
    arrow/.style={-Stealth, thick},
    tall_block/.style={base, minimum width=1.1cm, minimum height=1.8cm},
    std_block/.style={base, minimum width=0.9cm, minimum height=1.3cm},
    res_node/.style={base, minimum width=0.8cm, minimum height=1.1cm, font=\scriptsize, fill=white},
    plus/.style={circle, draw, thick, inner sep=0pt, minimum size=0.4cm, fill=white},
    enc_node/.style={base, fill=white, draw=blue!60!black, font=\scriptsize, minimum height=0.75cm},
    pool_node/.style={base, fill=orange!25, draw=orange!80!black, minimum width=2.2cm, font=\small\bfseries},
    gp_node/.style={base, fill=purple!15, draw=purple!80!black, minimum width=3.2cm, minimum height=1.1cm, font=\small},
    out_node/.style={base, fill=white, minimum width=5.5cm, minimum height=1cm, font=\scriptsize},
    info_box/.style={draw, dashed, inner sep=6pt, align=left, font=\tiny, fill=white, fill opacity=0.9}
}

\usepackage{algorithm}
\usepackage{algpseudocode}

\algnewcommand{\LineComment}[1]{\State\(\triangleright\) #1}

\usepackage{longtable}
\usepackage{tabularx}
\usepackage{hyperref}
\hypersetup{colorlinks=true, linkcolor=blue, citecolor=blue, urlcolor=red}

\begin{document}
\begin{frontmatter}
\title{A physics-informed neural network for improving surface reconstruction of intracranial saccular aneurysms via variational membrane equilibrium}

\author[cwnu]{Hyomin Ryu}

\author[scmc]{Seung Hwan Kim\corref{cor2}}
\ead{aajechtkskim@hanmail.net}

\author[cwnu]{Jaemin Kim\corref{cor1}}
\ead{jaeminkim@changwon.ac.kr}

\cortext[cor1]{Corresponding author}
\cortext[cor2]{Corresponding author}

\address[cwnu]{School of Mechanical Engineering, Changwon National University, Changwon 51140, Republic of Korea}
\address[scmc]{Department of Neurosurgery, Samsung Changwon Hospital, Sungkyunkwan University School of Medicine, Changwon 51353, Republic of Korea}
\journal{}

\begin{abstract}
Intracranial saccular aneurysms (ISAs) pose severe health risks, yet conventional population-based risk stratification scores (PHASES, UIATS, and ELAPSS) offer limited capacity for patient-specific rupture risk assessment. Image-based computational approaches have gained prominence, but traditional surface reconstruction relies on mathematical smoothing (e.g., L-curve criteria) that indiscriminately suppresses both imaging artifacts and genuine pathological features such as rupture-prone blebs. Although Laplace's membrane equilibrium ($\kappa_1 T_1 + \kappa_2 T_2 = P$) has long governed aneurysm wall mechanics (Humphrey and Kyriacou [Neurol. Res., 18 (1996)]), its integration into geometric reconstruction pipelines remains unexplored. This work introduces a physics-informed neural network (PINN) framework with B-spline representations, whose key contributions are: (i) embedding the variational equilibrium condition ($\delta \Pi = 0$) as a physics-informed loss that replaces the mathematical L-curve criterion with a biomechanically grounded artifact discrimination, (ii) developing a Manifold-Consistent CNN ansatz that preserves the closed-surface topology of vascular geometries, and (iii) establishing a Laplace equilibrium-driven reconstruction that filters imaging noise while preserving diagnostically critical high-curvature features. Application to patient-specific clinical datasets demonstrates that the framework eliminates non-physical concave artifacts without compromising genuine geometric anomalies. Clinical evaluation by a practicing neurosurgeon confirms that the resulting risk map---with rupture risk concentrated at the dome apex---is consistent with intraoperative observations, establishing a computational biomarker foundation for patient-specific rupture risk assessment of intracranial saccular aneurysms.
\end{abstract}


\end{frontmatter}

%%%%% Nomenclature Table %%%%%
\section*{Nomenclature}
\begingroup
\footnotesize
\setlength{\tabcolsep}{4pt}
\renewcommand{\arraystretch}{1.05}
\begin{center}
\begin{tabularx}{\textwidth}{l>{\raggedright\arraybackslash}X l>{\raggedright\arraybackslash}X}
\hline
Symbol & Definition & Symbol & Definition \\
\hline
$S(u,v)$ & NURBS surface mapping & $\mathbf{P}$, $\mathbf{P}_{macro}$, $\mathbf{P}_{base}$ & Control points (general, macro, upsampled) \\
$u, v$ & Parametric coordinates $\in [0,1]$ & $\mathbf{Q}$ & Measured point cloud \\
$N_{i,p}(u)$, $\mathbf{N}$ & B-spline basis functions / matrix & $\mathbf{d}$ & Control point displacement field \\
$U$, $V$ & Knot vectors & $w_{i,j}$ & NURBS weights (set to 1) \\
$\Pi$, $\delta \Pi$ & Total potential energy / first variation & $U_{stretch}$ & Membrane stretching energy \\
$W_{press}$ & External pressure work & $P$ & Transmural pressure \\
$\mathcal{L}_{total}$ & Total PINN loss function & $\mathcal{L}_{data}$ & Data fidelity loss \\
$\mathcal{L}_{equil}$ & Variational equilibrium loss & $\mathcal{L}_{tang}$ & Tangential suppression loss \\
$\lambda_{data}$, $\lambda_{phys}$, $\lambda_{tang}$ & Loss weighting factors & $\lambda_{reg}$, $\lambda_{rel}$ & Regularization / relative stiffness \\
$\alpha$, $\beta$ & Smoothing coefficients (tangential, curvature) & $\kappa_1$, $\kappa_2$ & Principal curvatures \\
$H$ & Mean curvature & $T_1$, $T_2$ & Principal wall tensions \\
$\sigma_1$, $\sigma_2$ & Principal stresses & $\sigma_R$ & Stress ratio (diagnostic marker) \\
$E$ & Young's modulus & $\nu$ & Poisson's ratio \\
$t$ & Wall thickness & $\hat{\mathbf{n}}$ & Unit normal vector \\
$\mathbf{g}_u$, $\mathbf{g}_v$ & Covariant basis vectors & $\mathbf{H}_E$ & Energy Hessian matrix \\
$\mathbf{C}$ & Kinematic constraint matrix & $\mathrm{SD}_{RMS}$ & RMS spatial deviation \\
$\mathcal{C}_{RMS}$ & RMS curvature error & $\Delta x_{img}$, $\Delta z_{img}$ & In-plane resolution / slice thickness \\
\hline
\end{tabularx}
\end{center}
\endgroup

\section{Introduction}

Intracranial saccular aneurysms (ISAs) represent a severe global health crisis; ruptured ISAs claim nearly 500,000 lives worldwide annually, with approximately 30,000 rupture events occurring in the United States alone \citep{kim2025vivo}. While population-based risk stratification scores (e.g., PHASES, UIATS, ELAPSS) are widely used, they offer limited capacity for patient-specific rupture risk assessment \citep{kim2025vivo, sturiale2021retrospective, pagiola2020phases}. Precise risk assessment necessitates accurately quantifying high-resolution localized wall stress and deformation patterns, as the roles of realistic geometry and material properties in governing aneurysmal mechanics have been extensively demonstrated~\citep{shah1997further, humphrey2000structure}. Consequently, reconstructing continuous, mathematically rigorous manifolds from clinical imaging is a fundamental prerequisite for reliable patient-specific biomechanical modeling.
In clinical workflows, this pipeline begins with Computed Tomography Angiography (CTA) or Magnetic Resonance Angiography (MRA). Raw 3D images are segmented using software like ScanIP~\citep{simpleware2024scanip} or 3D Slicer~\citep{fedorov2012slicer} to yield a preliminary discrete mesh. To leverage the ordered, centerline-based acquisition geometry of blood vessels rather than treating the data as an unorganized point cloud, vascular modeling typically extracts orthogonal cross-sections along the vessel centerline, allowing the preliminary mesh to be uniformly re-sampled layer-by-layer \citep{piccinelli2009framework}. This ordered extraction assumes that the vascular geometry is represented by a hierarchical sequence of structured point rings, yielding the structured ``point cloud data'' for subsequent parameterization.

However, the fidelity of this point cloud data is severely degraded by hardware-level limitations (e.g. angiographic machine), such as low signal-to-noise ratios (SNR) and partial volume effects. Paradoxically, angiography equipment post-processing--such as motion-compensation algorithms and reconstruction kernels--mathematically alters vascular boundaries, turning blurred boundaries into non-physical concave depressions \citep{berg2017does}. Because thin-walled membranes subject to internal pressure predominantly maintain convex geometries to sustain tensile stress, these concave depressions are mechanically invalid. When transferred to subsequent stress analyses, these geometric artifacts force the numerical system to generate spurious bending moments, triggering artificial stress concentration hotspots that corrupt rupture risk diagnostics \citep{valen2018realworld}. Historically, the computational mechanics community suppressed these geometric uncertainties using a three-tiered mathematical filtering framework: pixel-level anisotropic diffusion \citep{gerig1992nonlinear}, discrete mesh-level Laplacian smoothing \citep{taubin1995signal, cebral2010hemodynamics}, and control-point-level Tikhonov regularization optimized via L-curve criteria \citep{hansen1992analysis, styner2006framework}. Lacking physical constraints, these mathematical filters treat high-frequency artifacts and actual pathological features (such as rupture blebs \citep{ashkezari2021blebs, cebral2010bleb}) identically as high-curvature noise. Enforcing global smoothing penalties causes flattening sub-millimeter pathologies and triggering artificial volume shrinkage \citep{geng2013cortical, desbrun1999implicit}. As seen in \citet{cebral2010hemodynamics}, blindly obliterating blebs via discrete Laplacian filters induces spurious stress distributions, fundamentally compromising subsequent thermodynamic and hemodynamic predictions \citep{hoi2006validation}. To resolve these bottlenecks, smooth parametric surfaces—specifically Non-Uniform Rational B-Splines (NURBS)—have emerged for vascular modeling \citep{piegl1995nurbs, zhang2007patient}. Offering C² continuity, NURBS provides stable curvature fields essential for analytical stress formulations that depend heavily on high-order geometric derivatives. Unlike global polynomials \citep{farley2004parameterization}, the local support of NURBS represents localized blebs without global numerical oscillations. Although various pipelines automate NURBS fitting—including centerline sweeping \citep{morganti2015patient}, template registrations \citep{moola2026valvefit}, immersogeometric methods \citep{kamensky2015immersogeometric}, and learning-based autoencoders \citep{thani2026physics}—NURBS reconstruction specifically tailored for ISAs that preserves pathological landmarks while filtering noise remains unaddressed. 

Although Laplace's membrane equilibrium equation ($\kappa_1 T_1 + \kappa_2 T_2 = P$) has long been established as the fundamental governing relation for aneurysm wall mechanics~\citep{humphrey1996use}, its systematic integration into image-based geometric reconstruction via deep learning remains unexplored. To address this gap, the present work proposes an integrated workflow embedding physics-informed neural networks (PINN) \citep{karniadakis2021physics} directly into the geometric reconstruction pipeline. A robust reconstruction framework is introduced, synergizing a PINN \citep{raissi2019physics, saidaoui2024deep} with a Manifold-Consistent Convolutional Ansatz. By reformulating 3D surface reconstruction as a variational equilibrium problem governed by the principles of vascular mechanical homeostasis \citep{humphrey2000structure, baek2010clinical}, the PINN enforces physical filtering via variational equilibrium ($\delta \Pi = 0$). Instead of a mathematical L-curve, Laplace membrane equilibrium serves as the governing criterion to distinguish between non-physical artifacts (incapable of sustaining transmural pressure) and high-curvature features that maintain Laplace membrane equilibrium (true blebs). This physics-led strategy robustly preserves genuine diagnostic markers, improving simulation accuracy and consequently enhancing the predictive value of subsequent clinical-statistical indices.

To outline this approach, the remainder of this paper details the integrated computational framework. The methodology is organized as a sequential pipeline in Section 2. First, the overall architectural workflow and clinical data acquisition are established (Sections 2.1 and 2.2). The geometry reconstruction method is then formalized, detailing the constrained fitting for baseline template generation (Section 2.3). Building upon this, the manifold-consistent convolutional ansatz is introduced, where a neural network acts as a nonlinear trial function for ISA surface reconstruction. Finally, the framework is completed through a physics-informed loss formulation based on variational equilibrium ($\delta \Pi = 0$) to identify the optimal geometric configuration (Section 2.4). Section 3 presents the numerical validation.

\section{Theory and methodology}

\subsection{Overall framework pipeline}

\begin{figure*}[!t]
  \centering
  \includegraphics[width=\textwidth]{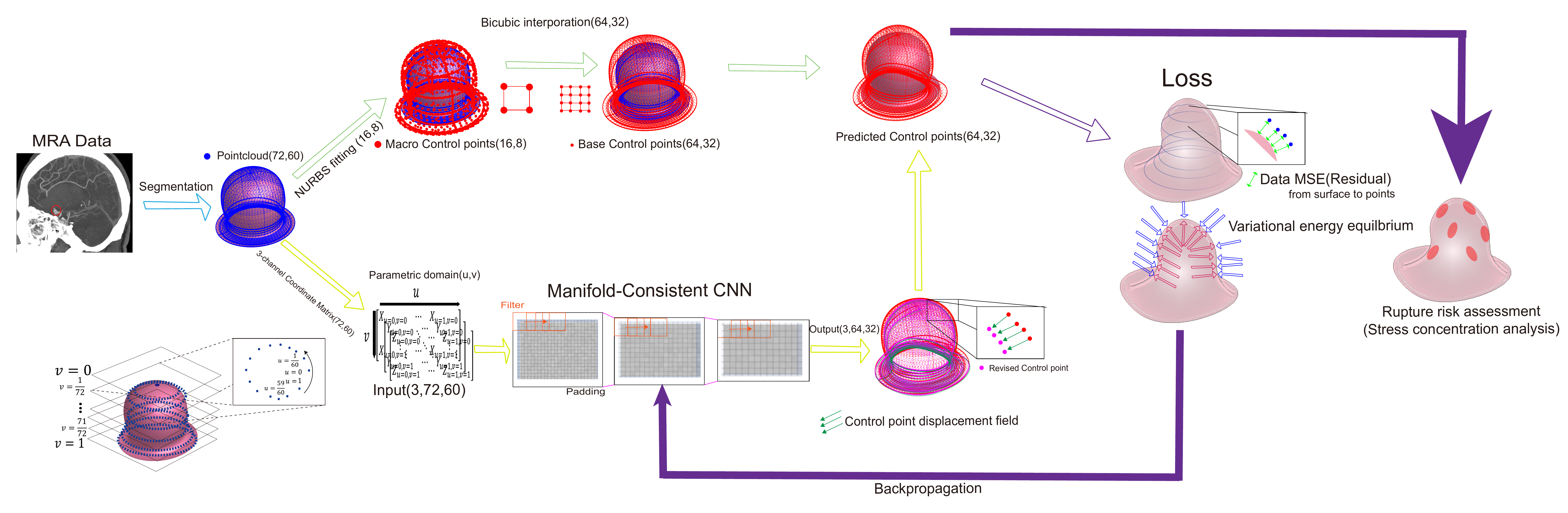}
  \caption{Schematic workflow of the proposed analysis-led modeling framework. The pipeline integrates clinical MRA data segmentation, Manifold-Consistent CNN prediction, and a physics-informed neural network governed by variational equilibrium ($\delta \Pi = 0$) to reconstruct artifact-free ISA geometries while preserving fine features.}
  \label{fig:overall_pipeline}
\end{figure*}

The proposed framework functions as a patient-specific numerical solver optimized to distinguish transient imaging artifacts from the genuine morphology of the ISA surface. By integrating B-spline parameterization with elastostatic solid mechanics, the framework acts as an active optimization process that drives the reconstructed geometry toward a strict Laplace membrane equilibrium state for each clinical case. Unlike conventional passive smoothing filters that compromise fine geometric details, this physics-driven approach ensures high fidelity in surface restoration. The overall structural workflow of the framework is depicted in Figure \ref{fig:overall_pipeline}.

As illustrated in the schematic pipeline, the procedure initiates with the raw point cloud data extracted via clinical MRA/CTA segmentation. This discrete geometric data is parameterized into a structured coordinate matrix within the parametric domain $(u, v)$ and concurrently utilized to extract macro control points via initial NURBS fitting. To construct a high-resolution geometric baseline, these macro control points are systematically upsampled through bicubic interpolation into base control points. Subsequently, the structured input tensor is processed through a specialized Manifold-Consistent CNN to predict the control point displacement field, yielding the updated and predicted control points. The central novelty lies in the active optimization loop driven by a dual-component loss function. Rather than relying on non-physical geometric metrics, the network updates its weights via backpropagation by continuously evaluating a data-fidelity mean squared error (MSE) residual alongside a variational equilibrium loss. This loss enforces the stationary condition of the total potential energy functional ($\delta \Pi = 0$) by explicitly balancing the thin-shell membrane stretching energy and the external workload driven by internal pressure. Through this integration of geometric parameterization and physical variational analysis, the architecture delivers a clean, artifact-free geometric profile optimized for highly reliable subsequent rupture risk assessment and stress multiaxiality analysis.

\subsection{Patient data acquisition and segmentation}

\begin{figure}[H]
  \centering
  \begin{subfigure}[b]{0.6\textwidth}
    \includegraphics[width=\linewidth]{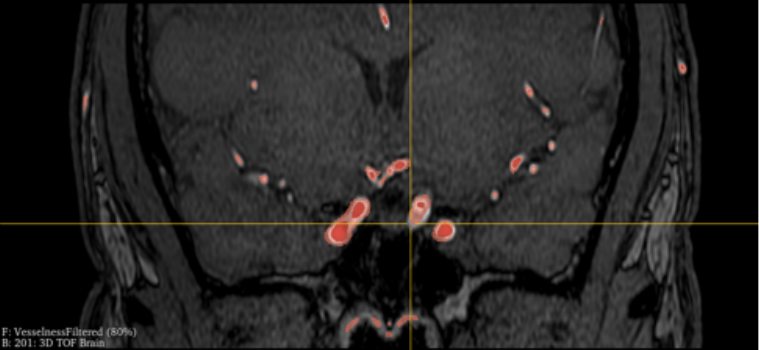}
    \caption{}
  \end{subfigure}
  \vspace{0.5cm}
  \begin{subfigure}[b]{0.6\textwidth}
    \includegraphics[width=\linewidth]{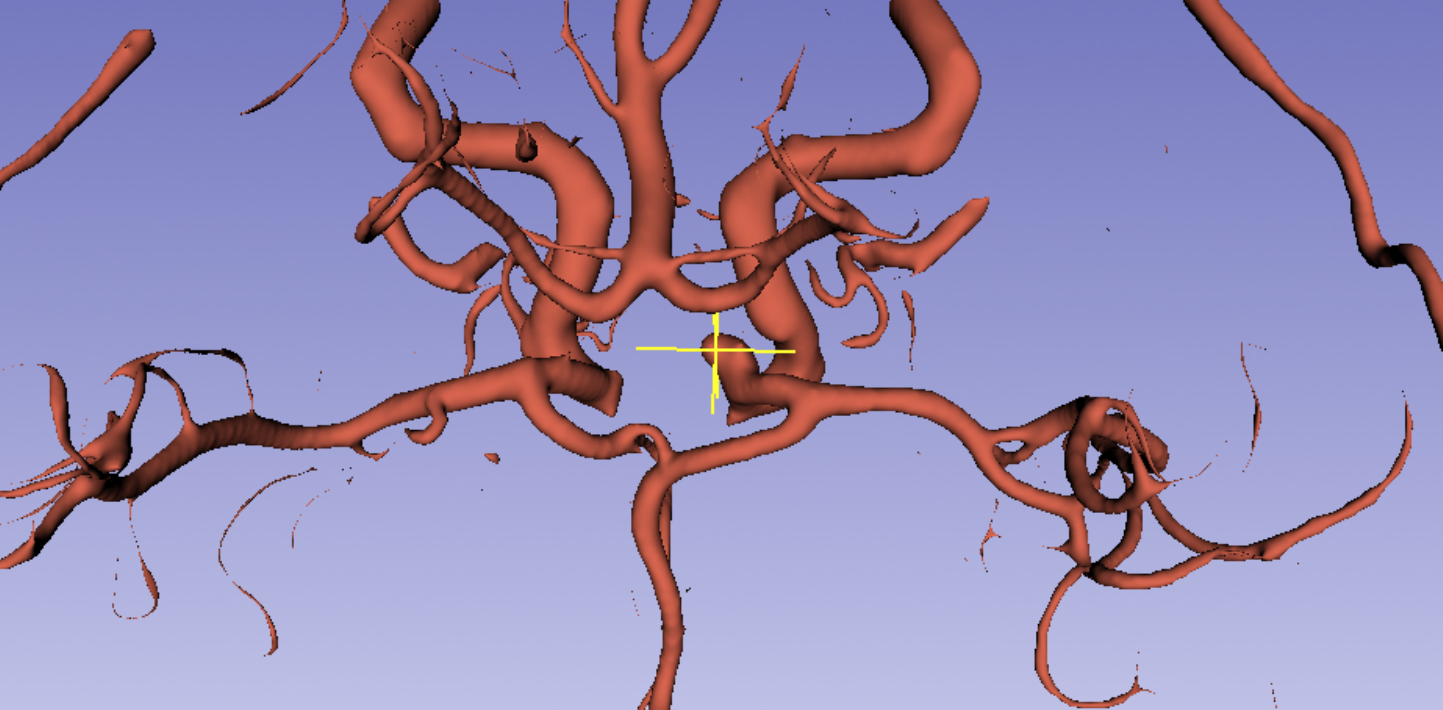}
    \caption{}
  \end{subfigure}
  \vspace{0.5cm}
  \begin{subfigure}[b]{0.6\textwidth}
    \includegraphics[width=\linewidth]{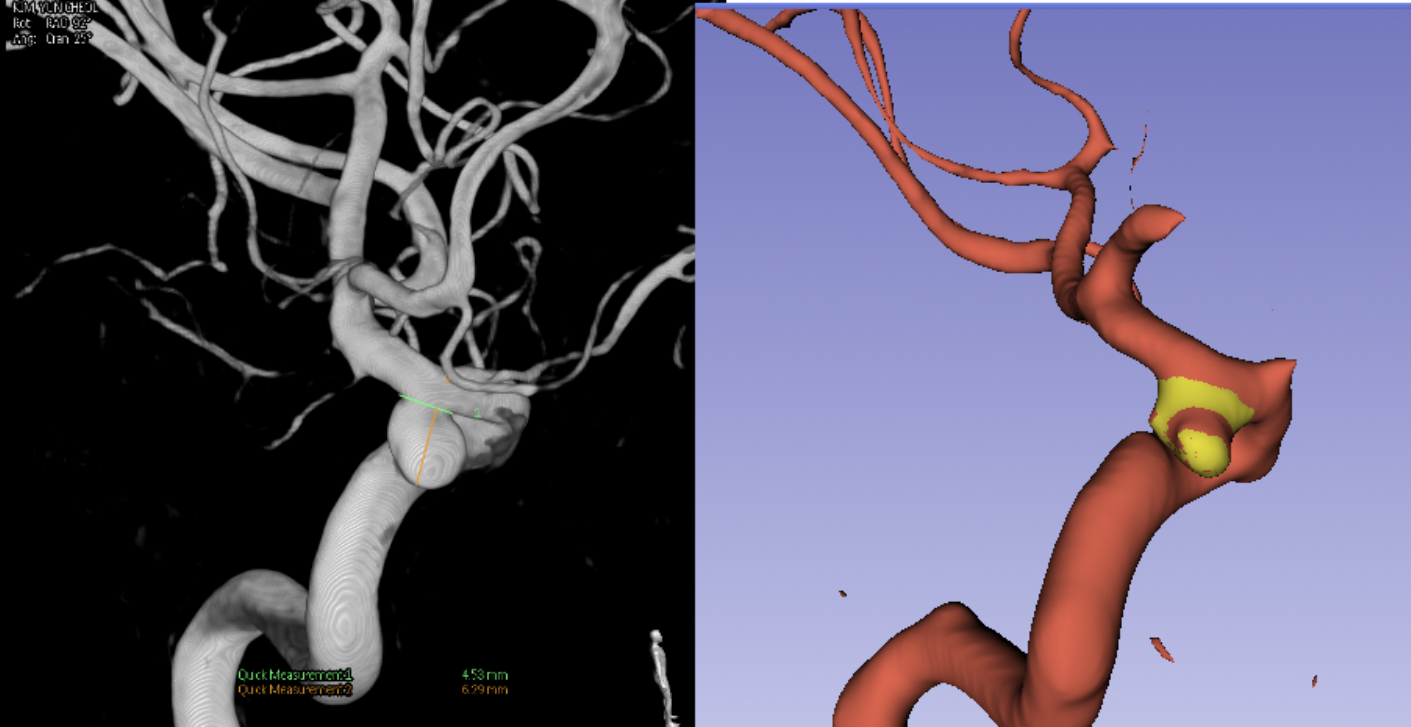}
    \caption{}
  \end{subfigure}
  \caption{Representative workflow of 3D geometric reconstruction in 3D Slicer. (a) Coronal slice of a 3D TOF Brain MRA with semi-automated vessel segmentation overlay. (b) 3D-rendered cerebrovascular anatomy (Circle of Willis) after segmentation; the yellow crosshair indicates the ISA location. (c) Left: clinical measurement of the aneurysm neck diameter ($D = 4.53$ mm) and height ($L = 6.29$ mm); Right: isolated ISA dome mesh (yellow) extracted using the Scissors tool.}
  \label{fig:slicer_workflow}
\end{figure}

\begin{figure}[H]
  \centering
  \includegraphics[width=0.98\linewidth]{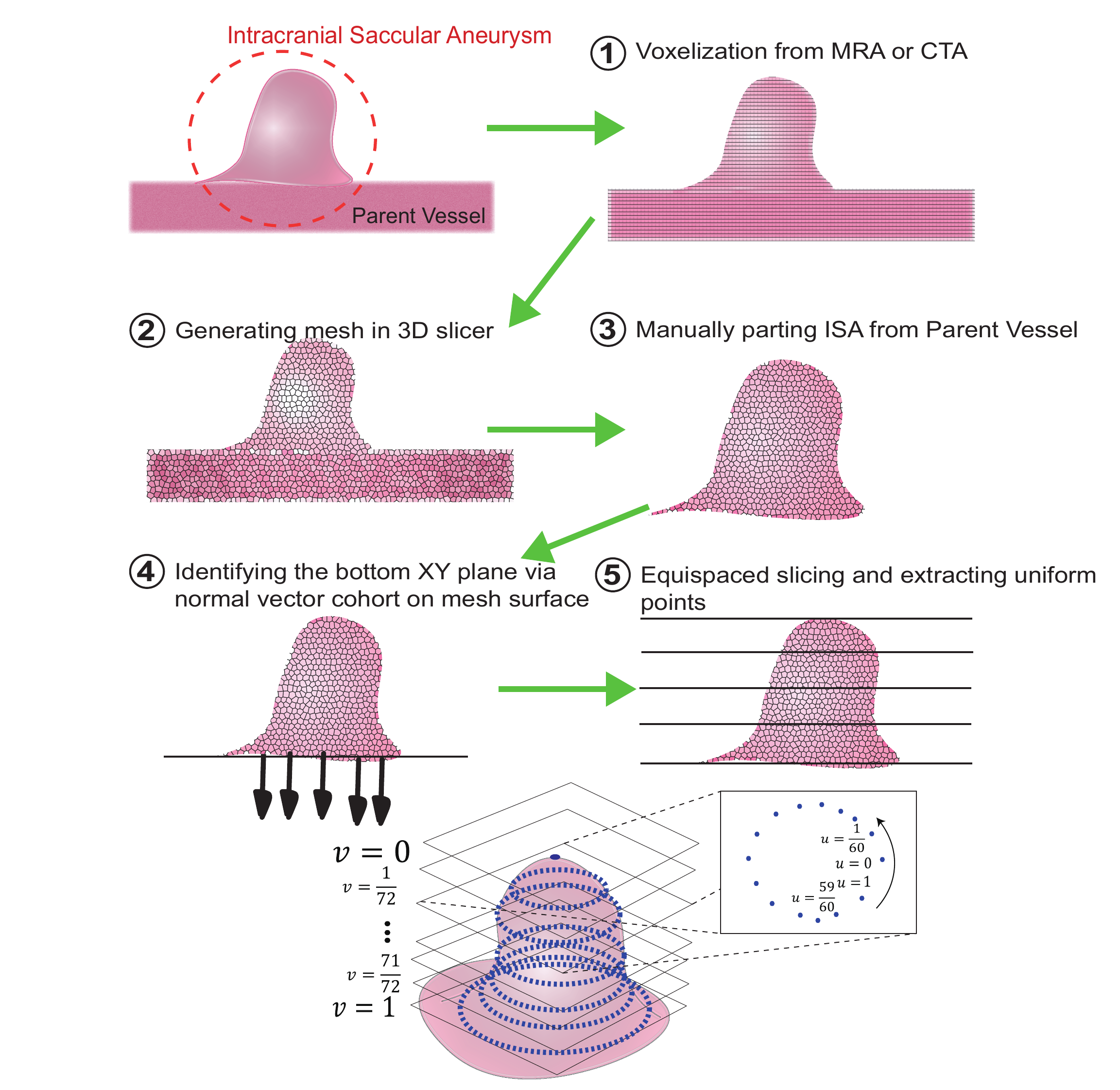}
  \caption{\small Schematic workflow of the unstructured point cloud extraction. \textcircled{1} Voxelization of the ISA and parent vessel from MRA or CTA imaging data. \textcircled{2} Surface mesh generation in 3D Slicer. \textcircled{3} Manual separation of the ISA dome from the parent vessel. \textcircled{4} Identification of the bottom (neck) plane via normal vector analysis on the mesh surface. \textcircled{5} Equispaced slicing along the meridional direction ($v = 0$ to $v = 1$) and extraction of uniformly distributed points along the circumferential direction ($u = 0$ to $u = 1$) on each slice.}
  \label{fig:extract_ptc}
\end{figure}

This study was conducted with the approval of the Institutional Review Board (IRB) of Samsung Changwon Hospital (IRB No. SCMC 2026-07-001), and written informed consent was obtained from all participants. The dataset explicitly included two distinct clinical scenarios based on standard high-resolution volumetric modalities: Computed Tomography Angiography (CTA) and Magnetic Resonance Angiography (MRA). CTA acquires data through X-ray attenuation by contrast agents, rendering it highly effective but prone to noise in a heavily corrupted Ruptured Cohort~\citep{hollingsworth2015reducing}. Conversely, non-contrast MRA leverages the flow-related enhancement of proton spins (e.g., Time-of-Flight), which intrinsically suffers from flow-induced signal voids and data sparsity in the Pre-Rupture Cohort~\citep{huston1994unruptured, bae2021compressed,lustig2008compressed}. Despite their distinct physical principles, both modalities fundamentally generate 3D voxel-based intensity maps stacked as sequential 2D slices within the standard RAS (Right-Anterior-Superior) coordinate system. These MRA and CTA volumetric datasets were then processed in 3D Slicer 5.10.0, utilizing semi-automated thresholding (0.07--0.3) and the `Scissors' tool to interpolate and segment the data, extracting the isolated preliminary discrete aneurysm dome mesh (\textbf{Figure~\ref{fig:slicer_workflow}}). To explicitly map these unstructured meshes into the spatial tensors required by the PINN architecture, an automated normal-driven alignment and grid sampling algorithm was implemented (\textbf{Figure~\ref{fig:extract_ptc}}). Unstructured point clouds were thereby mapped into a structured 3-channel spatial tensor ($X, Y, Z$) with dimensions of $3 \times 72 \times 60$.

\subsection{Geometry reconstruction}

The ISA surface is represented as a NURBS surface $S(u,v)$ over the parametric domain $u, v \in [0, 1]$~\citep{piegl1995nurbs, dierckx1995curve, bae2002nurbs}:

\begin{equation}
    S(u,v) = \frac{\sum_{i=0}^{n} \sum_{j=0}^{m} N_{i,p}(u) N_{j,q}(v) w_{i,j} \mathbf{P}_{i,j}}{\sum_{i=0}^{n} \sum_{j=0}^{m} N_{i,p}(u) N_{j,q}(v) w_{i,j}}
\end{equation}
where $N_{i,p}(u)$ and $N_{j,q}(v)$ are the B-spline basis functions of degree $p$ and $q$, defined over the knot vectors $U$ and $V$, respectively. The basis functions are recursively formulated via the Cox-de Boor algorithm~\citep{piegl1995nurbs, dierckx1995curve}:
\begin{equation}
    N_{i,0}(u) = \begin{cases} 1 & \text{if } u_i \leq u < u_{i+1} \\ 0 & \text{otherwise} \end{cases}
\end{equation}
\begin{equation}
    N_{i,p}(u) = \frac{u - u_i}{u_{i+p} - u_i} N_{i,p-1}(u) + \frac{u_{i+p+1} - u}{u_{i+p+1} - u_{i+1}} N_{i+1,p-1}(u) 
\end{equation} 

To parameterize the aneurysm surface properly, periodic and clamped uniform knot vectors are utilized for the azimuthal ($u$) and latitudinal ($v$) directions, respectively. Specifically, $U$ and $V$ are defined as a non-decreasing sequence of real numbers:

\begin{equation}
    U = \{\underbrace{u_0, \ldots, u_p}_{p+1}, u_{p+1}, \ldots, u_{n-p-1}, \underbrace{u_n, \ldots, u_{n+p}}_{p+1}\}
\label{eq:knot_vector_U} 
\end{equation}
\begin{equation}
    V = \{\underbrace{0, \ldots, 0}_{q+1}, v_{q+1}, \ldots, v_{m-q-1}, \underbrace{1, \ldots, 1}_{q+1}\}
\label{eq:knot_vector_V} 
\end{equation}

The term ``clamped'' implies that the first and last knot values are repeated $p+1$ times. For the ``periodic'' knot vector, the first and last $p+1$ knots are wrapped such that their corresponding knots fall outside the domain $[0, 1]$. This periodic wrapping guarantees $C^2$ continuity across the $u = 0 \rightarrow 1$ boundary, rendering $S(u,v)$ smoothly continuous along the circumferential direction. Additionally, both knot vectors are ``uniform,'' indicating that the interior knots are equally spaced.

\begin{figure}[H]
\centering
\begin{subfigure}[b]{0.49\textwidth}
    \centering
    \includegraphics[width=\linewidth]{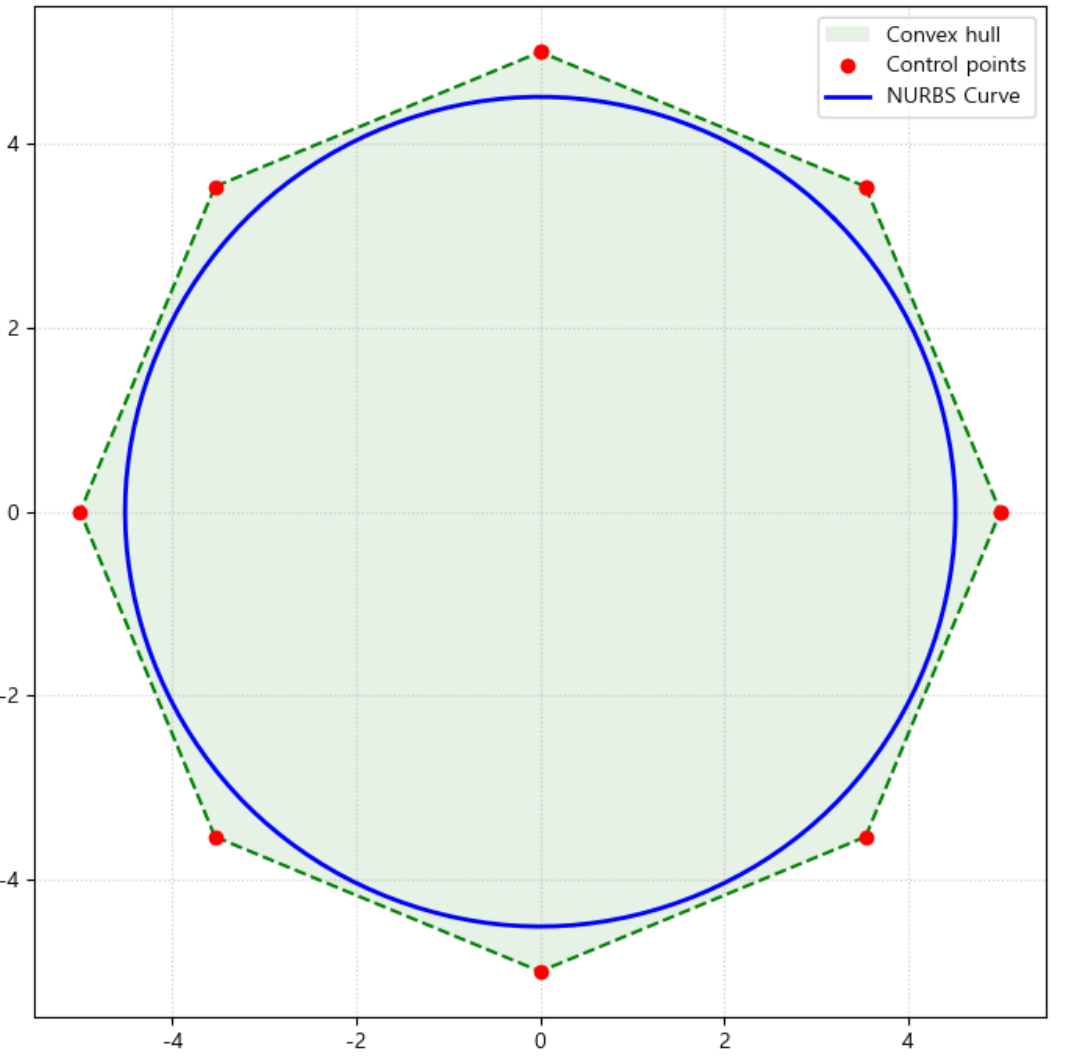}
    \caption{}
    \label{fig:side_by_side_1}
\end{subfigure}
\hfill
\begin{subfigure}[b]{0.49\textwidth}
    \centering
    \includegraphics[width=\linewidth]{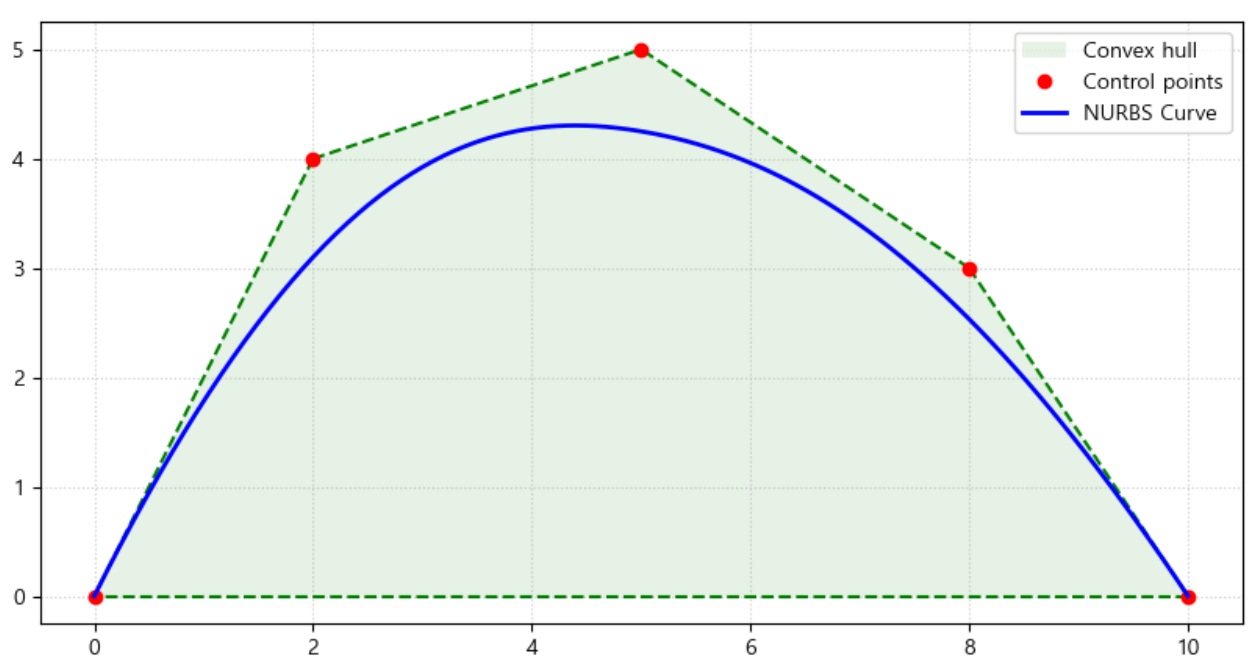}
    \caption{}
    \label{fig:side_by_side_2}
\end{subfigure}
\begin{subfigure}[t]{0.48\textwidth}
    \centering
    \includegraphics[width=\textwidth]{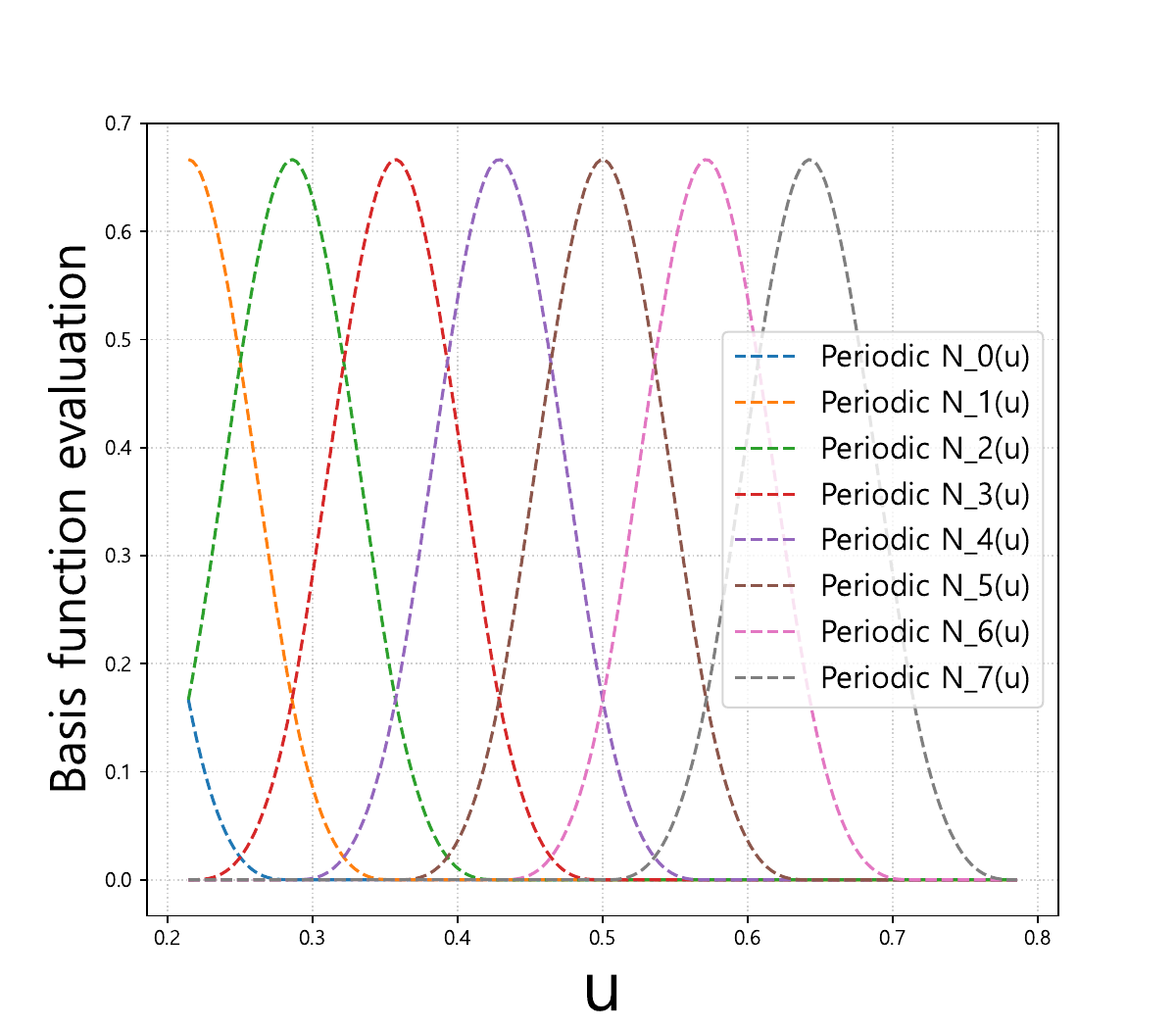}
    \caption{}
    \label{fig:sub_c}
\end{subfigure}
\hspace{\fill}
\begin{subfigure}[t]{0.48\textwidth}
    \centering
    \includegraphics[width=\textwidth]{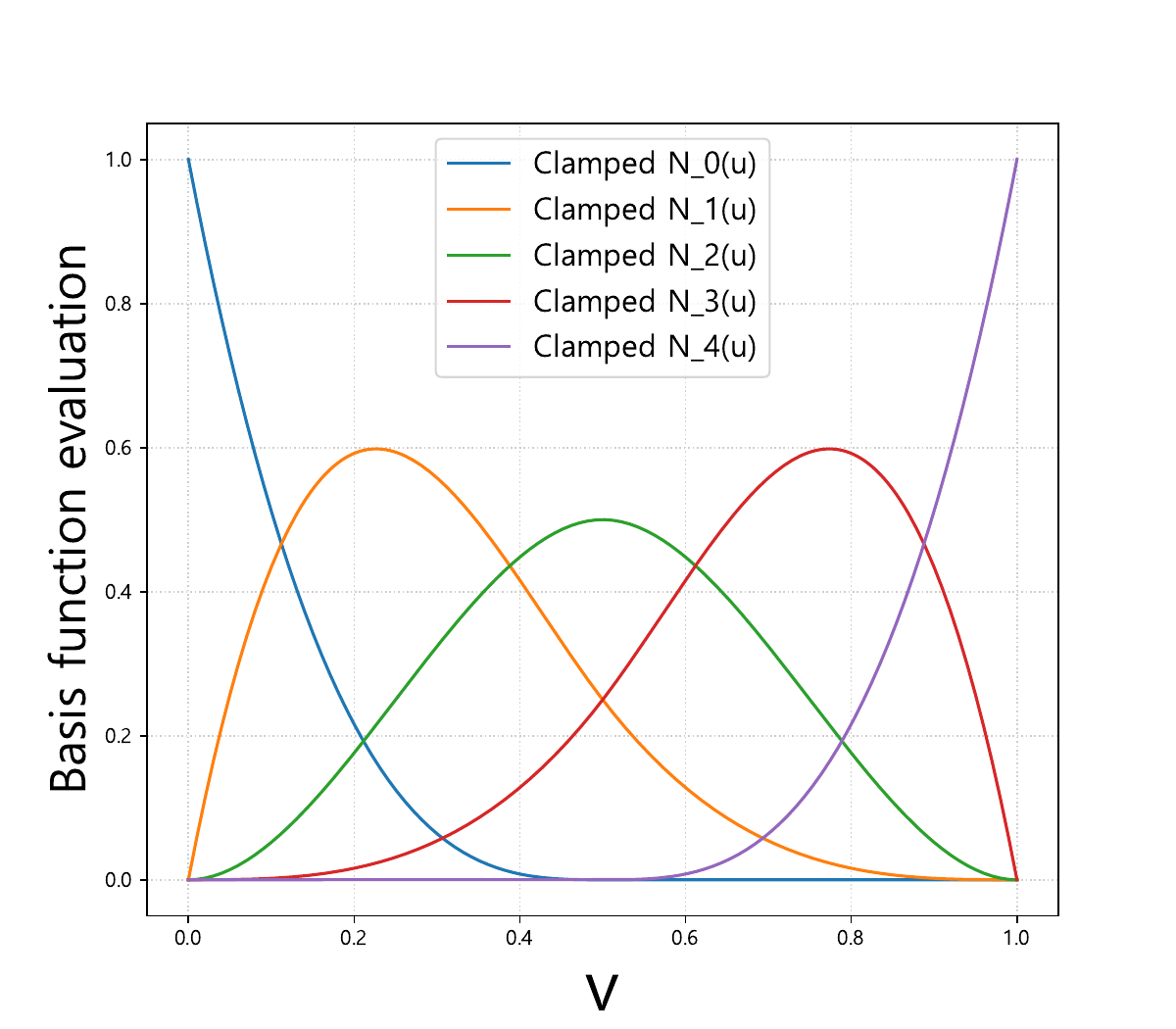}
    \caption{}
    \label{fig:sub_d}
\end{subfigure}
\caption{NURBS curve components. (a) Periodic curve with knot vector U \eqref{eq:knot_vector_U}, defining azimuthal direction in this paper. (b) Clamped curve with knot vector V \eqref{eq:knot_vector_V}, defining latitudinal direction in this paper. Evaluation of the B-spline basis functions that define the continuous surface. (c) In the azimuthal ($u$) direction, periodic knot vectors ensure $C^2$ continuity across the seam ($N(u=0)=N(u=1)$), (d) while in the latitudinal ($v$) direction, the basis functions are evaluated over an open knot vector to define the pole and boundary.}
\label{fig:side_by_side_example}
\end{figure}

In the proposed global geometry reconstruction, the homogeneous weight condition ($w_{i,j}=1$) is strictly enforced. This reduces the NURBS formulation to a standard non-rational B-spline surface~\citep{piegl1995nurbs, dierckx1995curve}, thereby reducing the fitting to a convex linear least-squares problem and preventing ill-conditioned gradients during subsequent neural backpropagation~\citep{raissi2019physics, wang2021understanding,tang2026neural}.

The fitting objective seeks to minimize the total energy functional $E(\mathbf{P})$, used here as a cost function, encompassing data fidelity, smoothing regularization ($\alpha, \beta$), and a Tikhonov regularization penalty ($\lambda_{reg}$) on the control points ($\mathbf{P}$)~\citep{lai2007unconstrained, hansen1992analysis}:

\begin{equation}\label{eq:fitting_objective}
    \min_{\mathbf{P}} E(\mathbf{P}) = \underbrace{\left\| \mathbf{Q} - \mathbf{N} \mathbf{P} \right\|^2}_{\text{Data fidelity}} + \underbrace{\alpha \left\| \mathbf{N}_\alpha \mathbf{P} \right\|^2}_{\text{Tangential smoothing}} + \underbrace{\beta \left\| \mathbf{N}_\beta \mathbf{P} \right\|^2}_{\text{Curvature smoothing}} + \underbrace{\lambda_{reg} \left\| \mathbf{P} \right\|^2}_{\text{Tikhonov reg.}}
\end{equation} 
where $\mathbf{N}$ denotes the matrix form of the B-spline basis functions $N_{i,p}(u)$, and $\mathbf{N}_\alpha$, $\mathbf{N}_\beta$ are the corresponding first- and second-derivative matrices with respect to the parameter $u$. The coefficients $\alpha$ and $\beta$ are the smoothing weights for the tangential and curvature terms, respectively, where $\beta$ is the key hyperparameter governing the smoothing--fidelity trade-off at the parametric fitting stage. Setting the optimality condition $\nabla_{\mathbf{P}} E = \mathbf{0}$ yields the normal equations with the Hessian matrix $\mathbf{H}_E$:

\begin{equation}\label{eq:hessian}
    \mathbf{H}_E = \mathbf{N}^T \mathbf{N} + \alpha \mathbf{N}_\alpha^T \mathbf{N}_\alpha + \beta \mathbf{N}_\beta^T \mathbf{N}_\beta + \lambda_{reg}\mathbf{I}
\end{equation} 

Unlike the smoothing terms ($\alpha$, $\beta$), which control geometric regularity, the Tikhonov regularization term ($\lambda_{reg} = 10^{-9}$) guarantees strict positive definiteness of $\mathbf{H}_E$ by shifting its eigenvalues, thereby preventing singularity in the matrix inversion~\citep{hansen1992analysis, bi2022span}. This fitting-based preprocessing operates at a resolution of $n_\theta \times n_\phi$.

The global geometry $\mathbf{P}_{macro} \in \mathbb{R}^{3 \times 16 \times 8}$ is resolved using the Lagrange multiplier method via the following block matrix inversion:

\begin{equation}\label{eq:block_matrix}
    \begin{bmatrix}
    \mathbf{H}_E & \mathbf{C}^T \\
    \mathbf{C} & \mathbf{0}
    \end{bmatrix}
    \begin{bmatrix}
    \mathbf{P}_{macro} \\
    \boldsymbol{\lambda}
    \end{bmatrix} =
    \begin{bmatrix}
    \mathbf{N}^T \mathbf{Q} \\
    \mathbf{D}
\end{bmatrix}
\end{equation}
where $\mathbf{C}$ is the constraint matrix enforcing pole convergence and $C^2$ continuity across the circumferential seam. By setting the smoothing weights to zero ($\alpha = \beta = 0$) at this stage, the global geometry reconstruction yields an unsmoothed reference configuration that preserves all measured features for the subsequent PINN-based fine-scale refinement. The resulting $\mathbf{P}_{macro}$ is then bicubically interpolated to $\mathbf{P}_{base} \in \mathbb{R}^{3\times64\times32}$ as shown in Figure~\ref{fig:overall_pipeline}.

This restriction to low-resolution fitting is motivated by computational cost: at high resolution (e.g., $64 \times 32$), the saddle-point system arising from the Lagrange multiplier formulation becomes prohibitively large to assemble and invert. Instead, the subsequent fine-scale refinement leverages PyTorch's automatic differentiation (Autograd) to enforce PDE constraints via gradient-based optimization, eliminating the need for explicit matrix assembly.

\begin{figure}[H]
    \centering
    \begin{subfigure}[t]{0.48\textwidth}
        \centering
        \includegraphics[width=\textwidth]{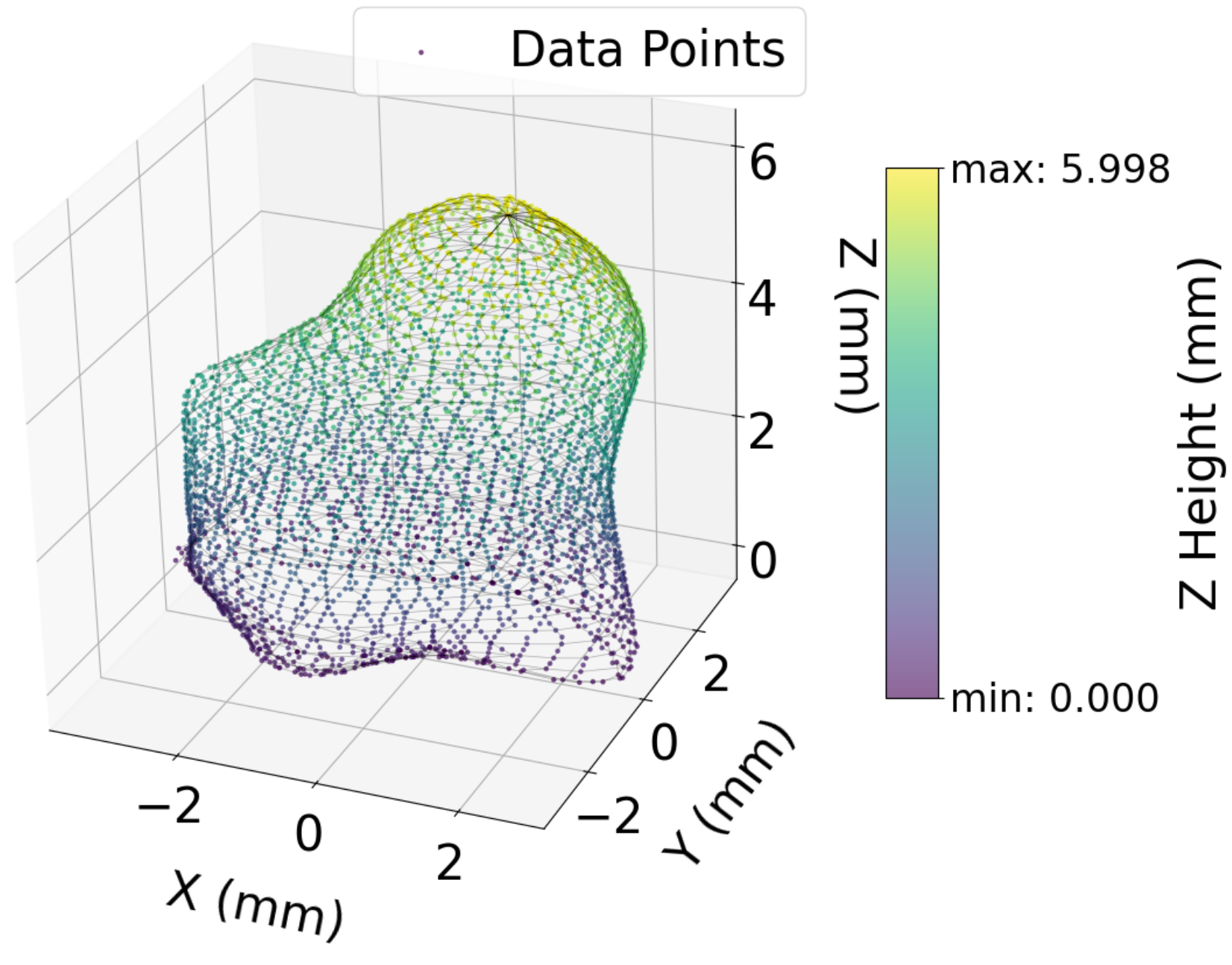}
        \caption{}
        \label{fig:sub_a}
    \end{subfigure}
    \hspace{\fill}
    \begin{subfigure}[t]{0.48\textwidth}
        \centering
        \includegraphics[width=\textwidth]{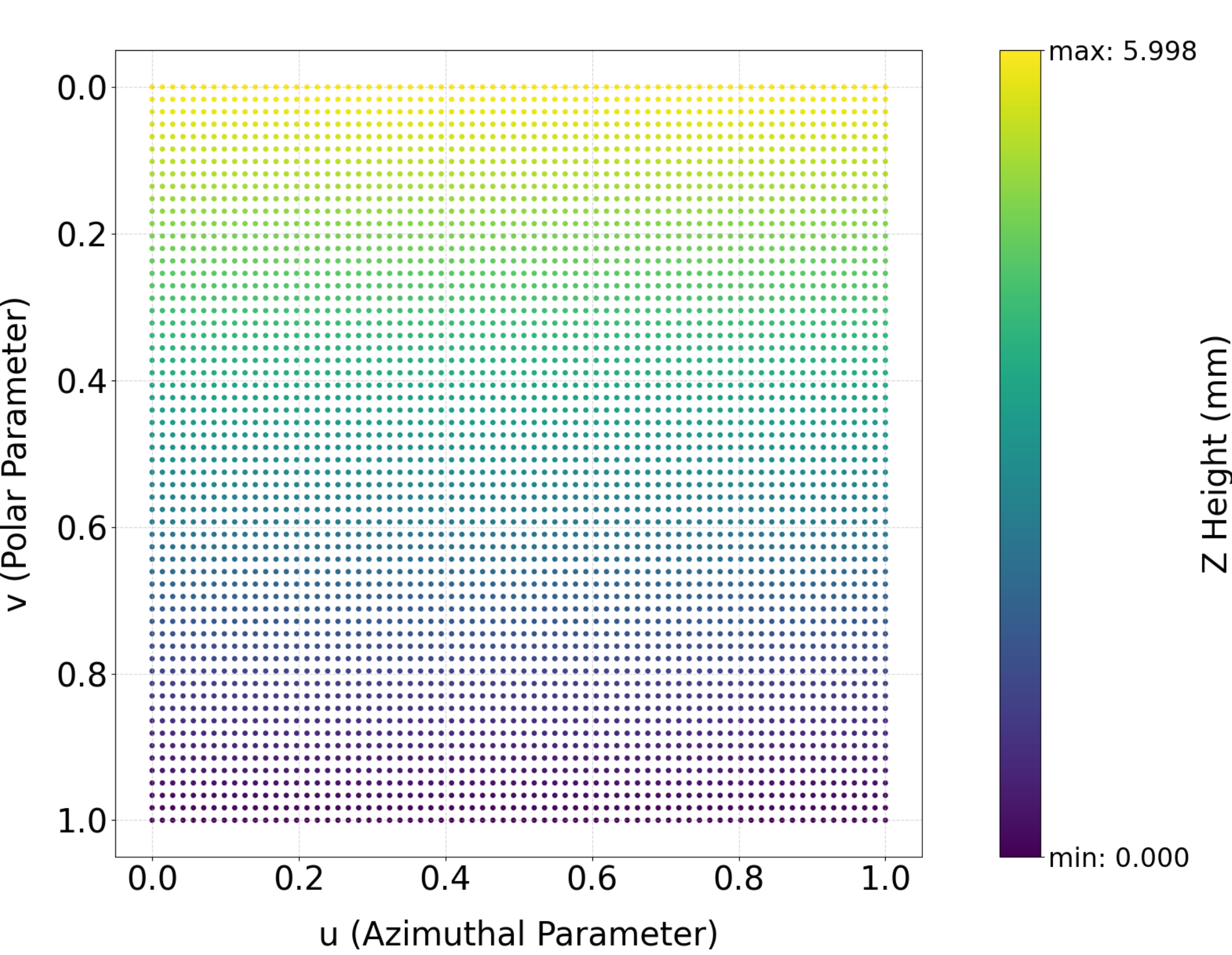}
        \caption{}
        \label{fig:sub_b}
    \end{subfigure}
    \caption{Components for surface reconstruction. (a) The raw ISA point cloud extracted from the mesh as described in Section 2.2. (b) The point cloud mapped to the parametric domain, illustrating the uniform sampling grid.}
    \label{fig:new_figure_group}
\end{figure}
\subsection{Physics-informed neural network architecture}

As discussed above, conventional smoothing trades off smoothness against residual errors, leading to an unavoidable loss of diagnostically critical high-curvature features. To address this limitation, the present study proposes a physics-informed neural network (PINN) framework~\citep{perdikaris2017nonlinear, kim2024review, hu2024physics} in which a Manifold-Consistent CNN serves as a nonlinear ansatz~\citep{haghighat2021physics, saidaoui2024deep}. As illustrated in Figure~\ref{fig:framework_structure}, the architecture progressively extracts spatial features ($3 \rightarrow 64$ channels) and predicts a control point displacement field ($64 \rightarrow 3$ channels), employing the $\tanh$ activation function for nonlinearity under a local axisymmetry assumption~\citep{humphrey1996use}.

\subsubsection{Manifold-consistent convolutional neural network}

The NURBS surface $S(u,v)$ is defined over a bi-parametric domain $(u, v) \in [0,1]^2$. The point extraction strategy aligns the segmented geometry along these parametric axes, producing a structured $72 \times 60$ coordinate matrix (Figure~\ref{fig:extract_ptc}). Conventional physics-informed solvers rely on fully connected multilayer perceptrons (MLPs)~\citep{raissi2019physics, saidaoui2024deep,haghighat2021physics, yang2021bpinns,hu2024physics,cai2021physicsfluid,cai2021physicsheat,misyris2020physics}, whose fully connected layers require flattening the input into a one-dimensional vector, inevitably destroying local spatial adjacency. To preserve these spatial correlations, a Convolutional Neural Network (CNN) architecture is employed~\citep{luo2016understanding, ronneberger2015unet}.

However, a na\"ive application of CNNs to this parametric grid introduces a fundamental inconsistency: the 2D grid is obtained by ``unfolding'' a closed 3D surface, so its edges are not true boundaries but artifacts of the parameterization. For instance, the left and right edges of the grid ($u=0$ and $u=1$) correspond to the same physical seam on the vessel wall, and the top edge ($v=0$) converges to a single pole at the dome apex. Standard CNN padding (zero, reflect, or circular) treats each edge uniformly, violating these heterogeneous topological constraints and producing spurious discontinuities in the predicted displacement field. The key idea of the proposed \textit{Manifold-Consistent CNN} is to embed these topological relationships directly into the padding operation, so that every convolutional kernel---even at the grid boundary---sees a neighborhood that is faithful to the original 3D manifold geometry.

\begin{figure}[H]
\centering
\includegraphics[width=\linewidth]{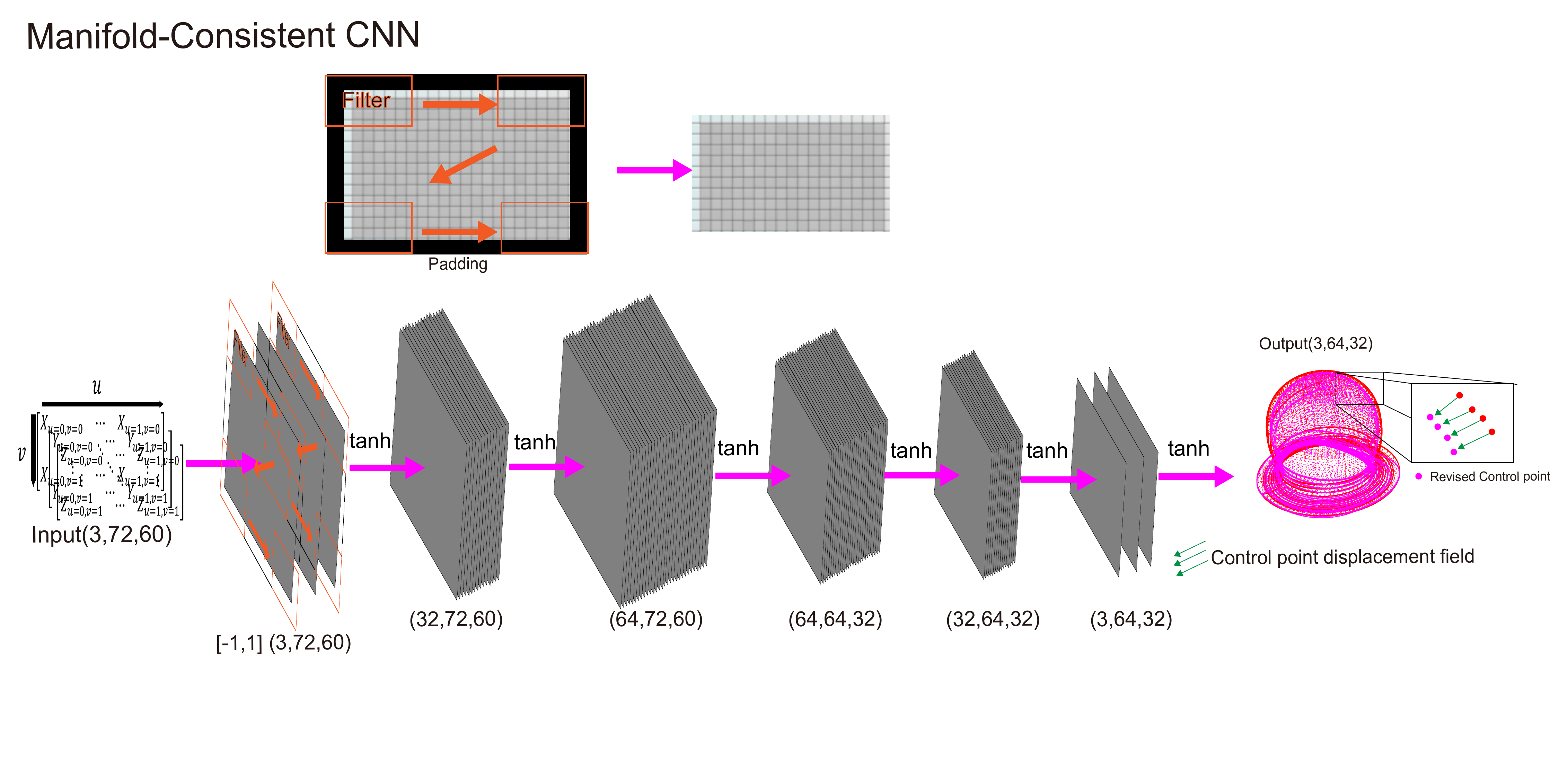}
\caption{Overview of the convolutional layers that extract spatial features (illustrated for the first channel) and predict the control point displacement field. The input coordinates are normalized to $[-1, 1]$ via max-absolute scaling to stabilize gradient propagation, and the final $\tanh$ activation is rescaled back to physical units.}
\label{fig:framework_structure}
\end{figure}

The proposed network (Figure~\ref{fig:framework_structure}) accepts the extracted 3-channel coordinate matrix ($3 \times 72 \times 60$) as input. Successive convolutional layers progressively downsample the spatial dimensions to $3 \times 64 \times 32$, matching the resolution of the target control points $\mathbf{P}$. The network parameters are optimized via backpropagation driven by a multi-objective loss function coupling data fidelity (point-to-surface residual) with a physics-informed equilibrium constraint.

Upon convergence, the network predicts the control point displacement field $\mathbf{d}$, which is added to the baseline control points $\mathbf{P}_{\text{base}}$:
\begin{equation}
\mathbf{P} = \mathbf{P}_{\text{base}} + \mathbf{d}
\end{equation}

\begin{figure}[H]
\centering
\includegraphics[width=\linewidth]{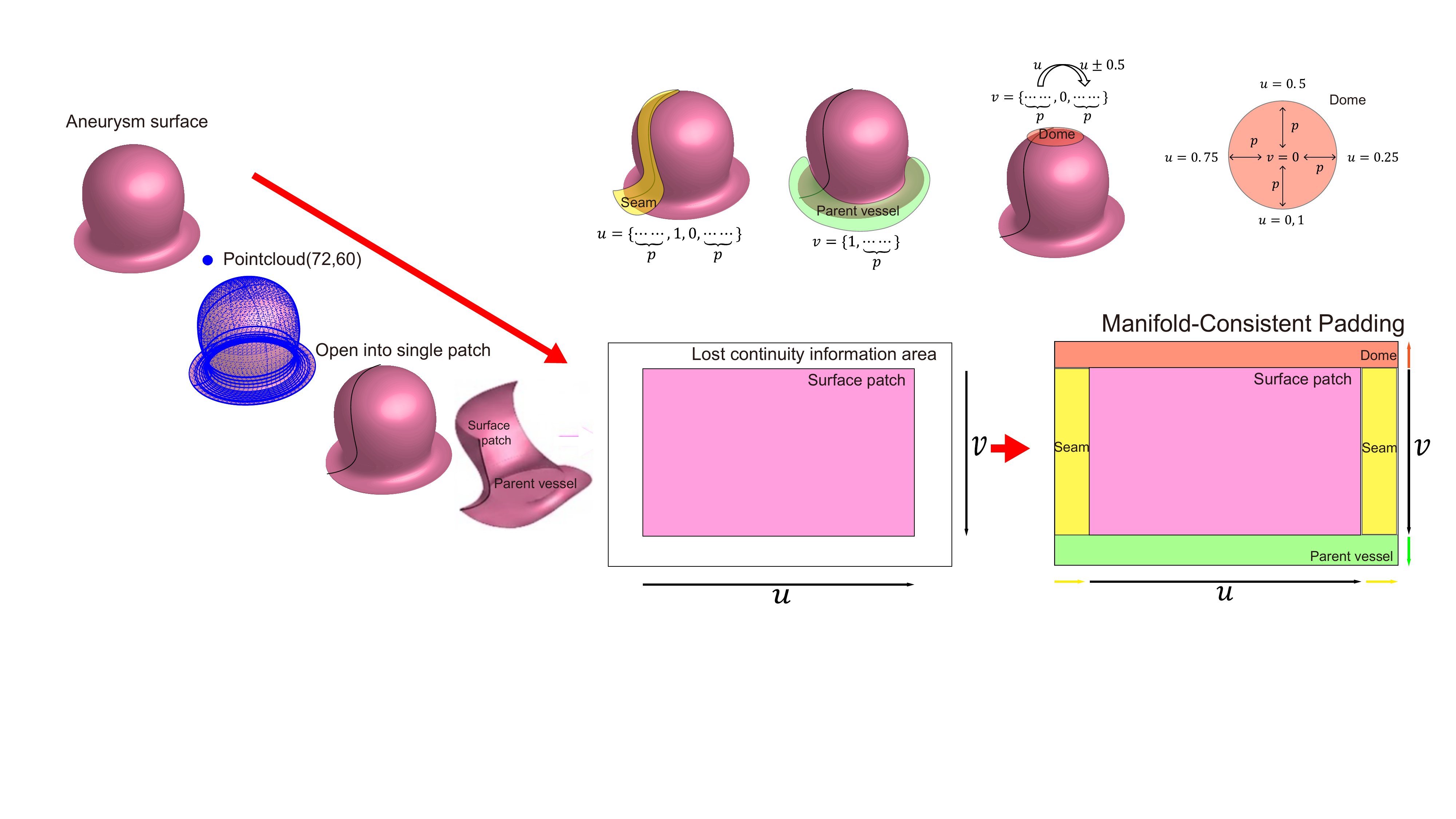}
\caption{Schematic of the Manifold-Consistent Padding. Flattening the closed 3D ISA manifold into a 2D patch disrupts periodic closure ($u$-axis), polar convergence ($v=0$), and parent vessel continuity ($v=1$). The specialized padding restores this connectivity with a width of $p+1$ knot spans, using interior value replication for the seam and dome boundaries and cubic extrapolation for the open neck boundary to suppress edge artifacts.}
\label{fig:backbone_detail}
\end{figure}

A major challenge in applying CNNs to closed manifolds is boundary padding. Standard padding schemes (e.g., zero, periodic, or symmetric) cannot accommodate the mixed boundary topology of the vascular grid: periodicity along the $u$-direction (seam closure), singular pole convergence at the dome apex ($v=0$), and open truncation at the neck ($v=1$). To handle these coexisting boundary types without inducing numerical artifacts, a custom Manifold-Consistent Padding scheme is developed (Figure~\ref{fig:backbone_detail}).

Mapping the 3D closed surface onto a 2D parametric patch disconnects intrinsic spatial relationships. The proposed padding scheme restores this continuity by imposing boundary conditions consistent with the NURBS geometric continuity. The padding width is set to $p+1$ to match the local support of the B-spline basis functions of degree $p$. Specifically, the periodic seam ($u=0,1$) and polar convergence ($v=0$) are handled by replicating adjacent interior values, while the open neck boundary ($v=1$) employs cubic extrapolation to suppress boundary artifacts.

\subsubsection{Loss function}
The PINN framework determines the optimal surface geometry by minimizing the total loss $\mathcal{L}_{total}$~\citep{raissi2019physics, cai2021physicsheat, yang2021bpinns}, which balances data fidelity, physics-informed equilibrium, and tangential displacement suppression:
\begin{equation}
    \mathcal{L}_{total} = \underbrace{\lambda_{data} \mathcal{L}_{data}}_{\text{Data fidelity}} + \underbrace{\lambda_{phys} \mathcal{L}_{equil}}_{\text{Equilibrium constraint}} + \underbrace{\lambda_{tang} \mathcal{L}_{tang}}_{\text{Tangential suppression}}
\end{equation}

\paragraph{(1) Data fidelity loss ($\mathcal{L}_{data}$)}
This term enforces agreement between the reconstructed NURBS surface $S(u,v)$ and the measured point cloud $\mathbf{Q}$~\citep{gao2021phygeonet,moola2026valvefit}. By penalizing the point-to-surface distance in physical space, it constrains the solution space to geometries consistent with the clinical measurements:
\begin{equation}\label{eq:data_loss}
    \mathcal{L}_{data} = \sum_{k=1}^{\rho_\theta} \sum_{l=1}^{\rho_z} \left\| S(u_k, v_l) - \mathbf{Q}_{k,l} \right\|^2
\end{equation}

\paragraph{(2) Variational equilibrium loss ($\mathcal{L}_{equil}$)}
ISA formation involves progressive deformation of the vessel wall under transmural pressure. The stationarity condition of the total potential energy ($\delta \Pi = 0$), which governs this mechanical equilibrium, is incorporated as a physics-informed loss to reconstruct a physically admissible shape while suppressing imaging artifacts~\citep{raissi2019physics, haghighat2021physics}. Under the membrane assumption, the loss is defined as the squared gradient norm of the total potential energy with respect to the displacement field, enforcing the balance between internal stretching ($U_{stretch}(\mathbf{d})$) and external pressure work ($W_{press}(\mathbf{d})$):
\begin{equation}\label{eq:equil_loss}
    \mathcal{L}_{equil} = \left\| \nabla_{\mathbf{d}} \left( \lambda_{rel} U_{stretch}(\mathbf{d}) - W_{press}(\mathbf{d}) \right) \right\|^2
\end{equation}
where $U_{stretch}(\mathbf{d}) \approx \iint \|\nabla \mathbf{d}\|^2 \, dA$ and $W_{press}(\mathbf{d}) \propto \iint (\kappa_1(\mathbf{d}) + \kappa_2(\mathbf{d})) \, dS$. The relative stiffness coefficient $\lambda_{rel}$ ensures dimensional consistency between the internal restoring force and the external pressure load. The detailed formulation is provided in the Appendix.

\paragraph{(3) Tangential displacement suppression ($\mathcal{L}_{tang}$)}
The tangential displacement loss $\mathcal{L}_{tang}$ penalizes control point sliding along the tangent plane, constraining the displacement to evolve strictly in the normal direction. This restriction is grounded in the Hadamard--Zol\'{e}sio structure theorem of shape sensitivity analysis, which states that only the normal component of a boundary displacement contributes to a physical change in the domain geometry, while tangential variations merely reparameterize the surface~\citep{sokolowski1992introduction}:
\begin{equation}\label{eq:tangent_loss}
    \mathcal{L}_{tang} = \sum_{i=1}^{64\times32} \left\| \mathbf{d}_{tangent}^{(i)} \right\|^2
\end{equation}

These three loss terms play complementary roles. The data fidelity loss $\mathcal{L}_{data}$ anchors the reconstructed surface to the patient-specific point cloud, preventing the physics loss from driving the geometry away from the clinical measurements. The equilibrium loss $\mathcal{L}_{equil}$ enforces the stationarity of total potential energy ($\delta \Pi = 0$), filtering out concave artifacts that cannot sustain transmural pressure while preserving genuine high-curvature features (e.g., blebs) that satisfy Laplace membrane equilibrium. The tangential suppression loss $\mathcal{L}_{tang}$ ensures that every control point update produces a physical shape change rather than a mere reparameterization of the surface.

This combination addresses a limitation that conventional smoothing cannot resolve: mathematical filters (e.g., L-curve criteria) treat all high-curvature regions uniformly as noise, whereas the physics-informed loss provides a biomechanically grounded discriminator---only geometries incompatible with the mechanical equilibrium are suppressed. Starting from the template $\mathbf{P}_{base}$, the iterative optimization thus emulates the growth-induced deformation of the aneurysm wall under transmural pressure. The underlying formulation adopts a set of biomechanical assumptions~\citep{sadasivan2013physical, watton2009modelling, yahya2010three} (detailed in Appendix B) that reduce the complex nonlinear vascular deformation to a tractable variational equilibrium problem.

\section{Numerical examples and validation}

The proposed PINN framework is validated on a representative clinical case of an unruptured ISA. Section~\ref{sec:risk_metrics} defines the biomechanical evaluation metrics---the stress ratio ($\sigma_R$) and mean curvature ($H$). Section~\ref{sec:artifact_removal} demonstrates the progressive removal of non-physical concave artifacts (Figure~\ref{fig:curvature_change}) and the resulting improvement in the rupture risk map (Figure~\ref{fig:risk_progression}). Section~\ref{sec:comparison} compares the PINN-optimized surface against conventional L-curve-based smoothing (Figure~\ref{fig:L-curve}) and evaluates the diagnostic reliability of the resulting risk maps (Figure~\ref{fig:risk_comparison}).

\subsection{Biomechanical risk assessment and evaluation metrics}\label{sec:risk_metrics}

The diagnostic utility of the reconstructed surfaces was evaluated using the stress ratio ($\sigma_R$). According to the Laplace membrane equilibrium~\citep{humphrey1996use, humphrey2000structure, baek2005competition, kim2025phasefield}, the principal wall tensions ($T_1, T_2$) are related to the principal curvatures ($k_1, k_2$) and the transmural pressure ($P$) as:
\begin{equation}
    T_1 = \frac{P}{2 \kappa_2}, \quad T_2 = \frac{P}{\kappa_2} \left( 1 - \frac{\kappa_1}{2 \kappa_2} \right)
\end{equation}
For an isotropic membrane with locally uniform thickness $t$, the principal stresses are $\sigma_i = T_i / t$. Because the stress ratio is defined as $\sigma_2 / \sigma_1$, the thickness $t$ and pressure $P$ cancel out:
\begin{equation}
    \sigma_R(u,v) = \frac{\sigma_2}{\sigma_1} = \frac{T_2}{T_1} = 2 - \frac{\kappa_1}{\kappa_2}
\end{equation}
Consequently, $\sigma_R$ is a purely geometry-derived indicator of stress multiaxiality, independent of material properties and loading magnitude. For visualization, this study introduces a logarithmic min-max normalization to map $\sigma_R$ onto a $[0, 100]$ percentage scale (risk score):
\begin{equation}
    \text{Risk Score}(u,v) = \frac{\log(\sigma_R(u,v)) - \log(\sigma_{R,\min})}{\log(\sigma_{R,\max}) - \log(\sigma_{R,\min})} \times 100
\end{equation}

In addition, the mean curvature ($H$) was used as a complementary metric to verify the morphological validity of the reconstructed surface, ensuring the elimination of non-physical concave artifacts ($H < 0$):
\begin{equation}
    H(u,v) = \frac{\kappa_1(u,v) + \kappa_2(u,v)}{2}
\end{equation}

\subsection{Artifacts removal and risk map optimization}\label{sec:artifact_removal}
Although the Neo-Hookean model does not fully capture the complex mechanical behavior of the arterial wall~\citep{humphrey2000structure, humphrey1996use, humphrey2008intracranial}, it is sufficient as a foundational constitutive model to establish a key physical constraint: a thin-walled membrane under positive internal pressure---analogous to a rubber balloon~\citep{boo2026isogeometric}---cannot sustain concave surface profiles under Laplace membrane equilibrium. Because ISA pathogenesis shares this loading condition, the mean curvature must remain strictly positive ($H > 0$) on the aneurysm wall.

\begin{figure}[H]
    \centering
    \includegraphics[width=1\linewidth, height=0.28\linewidth]{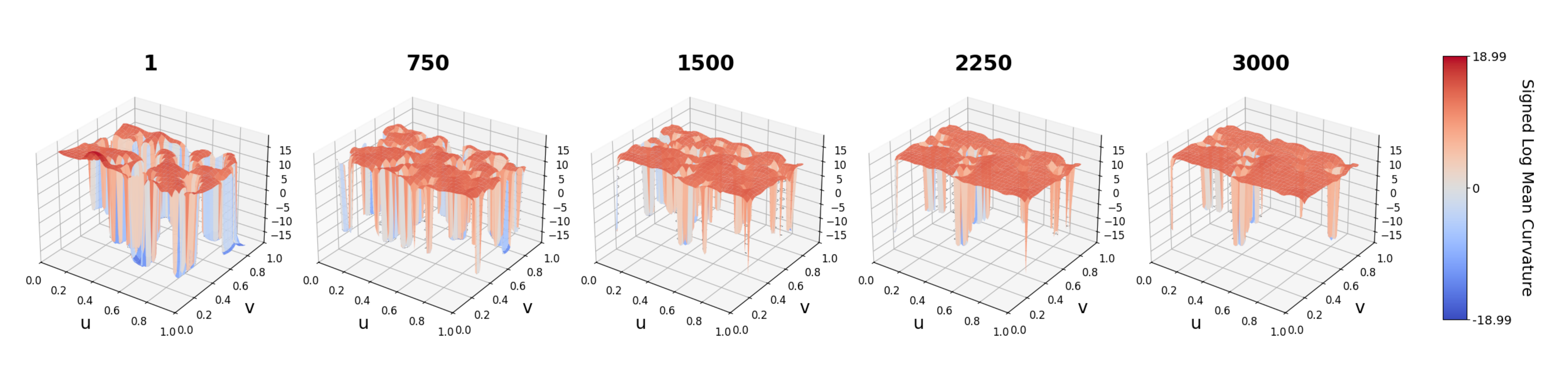}    \caption{Evolution of signed log mean curvature profiles over the parametric domain $(u, v)$ at iterations 1, 750, 1500, 2250, and 3000. At iteration 1, the baseline surface exhibits numerous downward spikes (blue, $H < 0$) corresponding to non-physical concave artifacts inherited from the constrained fitting. As the PINN optimization progresses, the equilibrium loss ($\mathcal{L}_{equil}$) progressively suppresses these negative curvature regions, while the positive curvature landscape (red/orange) representing the convex dome geometry is preserved. By iteration 3000, virtually all negative spikes are eliminated, yielding a curvature field consistent with a pressurized membrane ($H > 0$).}
    \label{fig:curvature_change}
\end{figure}

\begin{figure}[H]
    \centering
    \includegraphics[width=1\linewidth, height=0.28\linewidth]{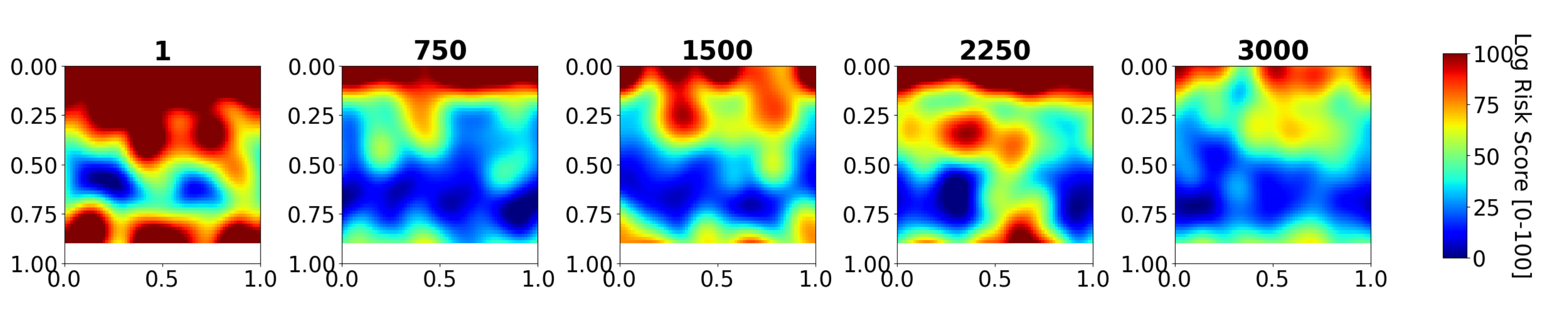}    
    \caption{Evolution of the log risk score maps over the parametric domain $(u, v)$ at iterations 1, 750, 1500, 2250, and 3000. At iteration 1, the risk map is dominated by saturated high-risk zones (dark red) near the dome apex ($v \approx 0$) caused by the concave artifacts shown in Figure~\ref{fig:curvature_change}. As the optimization proceeds, these artifact-induced hotspots are progressively resolved, and the risk distribution converges to a spatially coherent pattern. By iteration 3000, the high-risk concentration is localized at the dome apex region, consistent with clinical expectations for ISA rupture sites.}
    \label{fig:risk_progression}
\end{figure}

As shown in Figure~\ref{fig:curvature_change}, the initial surface derived from the constrained fitting baseline exhibits numerous concave artifacts ($H < 0$), which are non-physical for pressurized ISAs. These artifacts generate spurious stress concentration zones erroneously scattered across the vascular wall. Figure~\ref{fig:risk_progression} demonstrates the progressive correction through PINN optimization: by enforcing variational equilibrium ($\delta \Pi = 0$), the framework eliminates negative curvature regions while preserving genuine high-curvature features. The resulting risk map is free of artificially scattered hotspots, concentrating the risk at biomechanically meaningful locations.

\subsection{Comparison with conventional methods}\label{sec:comparison}

\begin{figure}[H]
    \centering
    \begin{subfigure}[t]{0.48\textwidth}
        \centering
        \includegraphics[width=\textwidth]{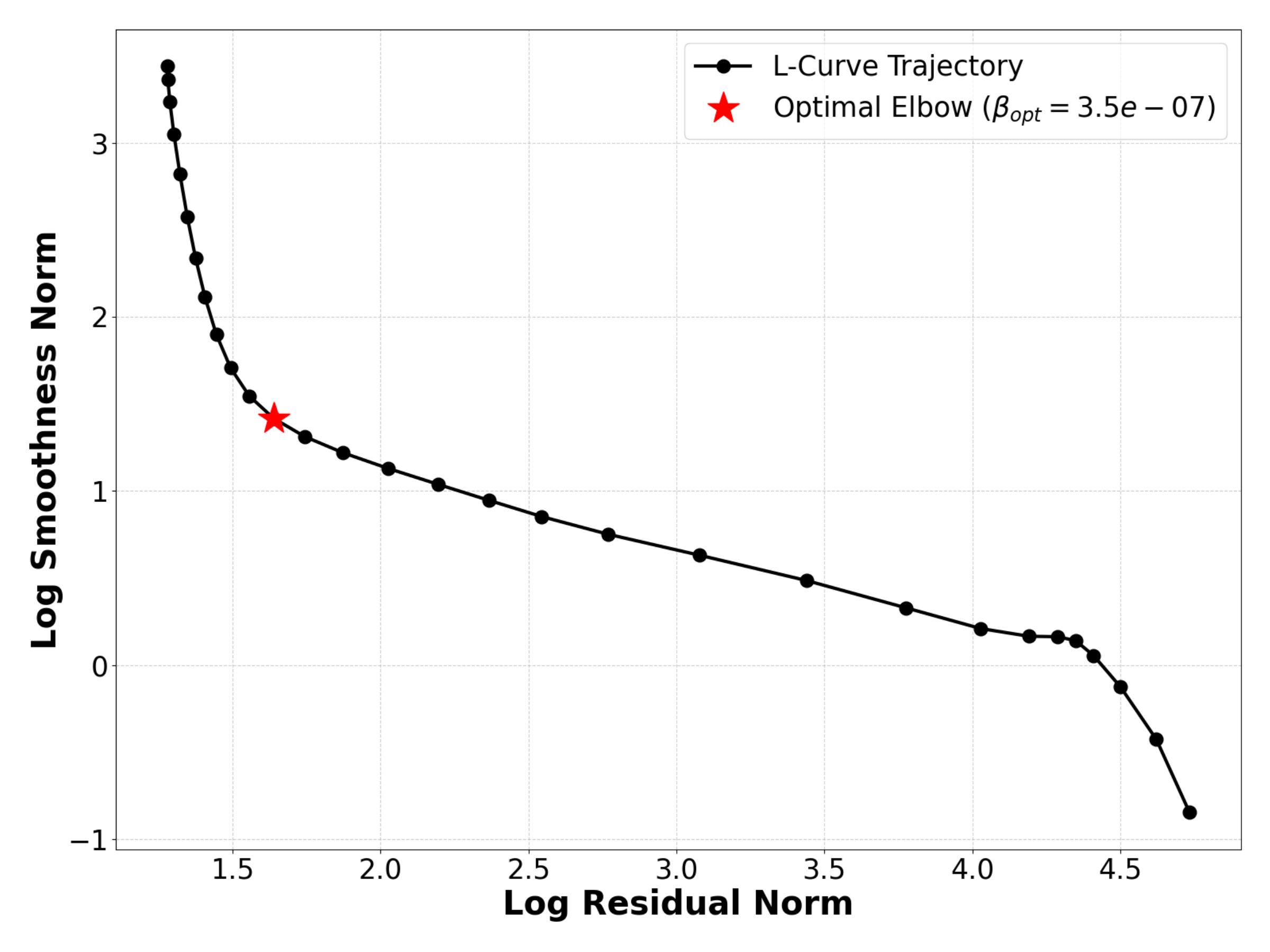}
        \caption{}
        \label{fig:sub_E}
    \end{subfigure}
    \hspace{\fill}
    \begin{subfigure}[t]{0.48\textwidth}
        \centering
        \includegraphics[width=\textwidth]{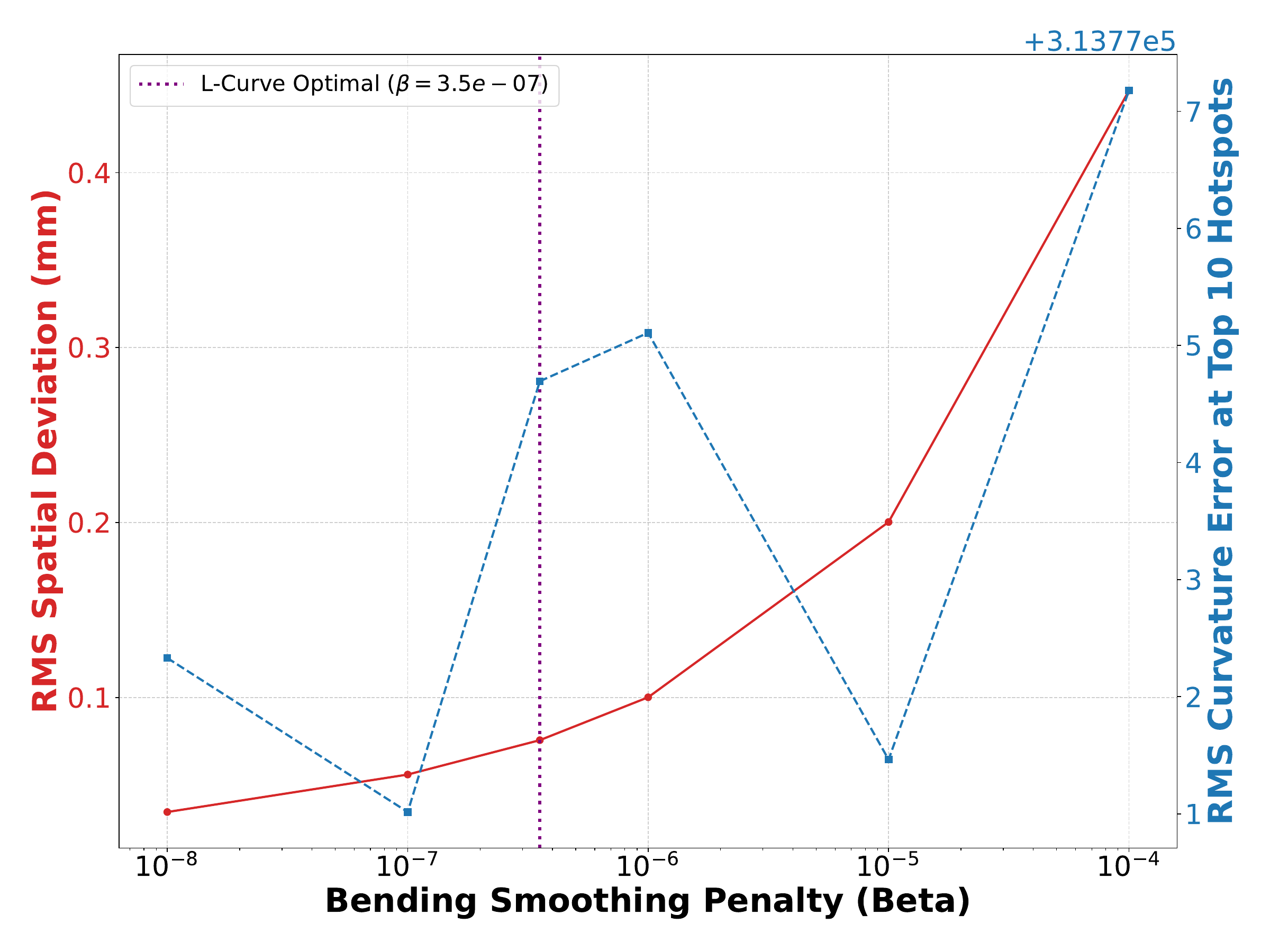}
        \caption{}
        \label{fig:sub_f}
    \end{subfigure}
    \caption{L-curve analysis for the smoothing factor $\beta$. (a)~L-curve trajectory in log--log space; the red star marks the optimal elbow at $\beta_{\mathrm{opt}} = 3.5 \times 10^{-7}$. (b)~$\mathrm{SD}_{\mathrm{RMS}}$ (red solid, left axis) and $C_{\mathrm{RMS}}$ (blue dashed, right axis) as functions of $\beta$; the vertical dotted line indicates $\beta_{\mathrm{opt}}$. The persistent offset of $+3.1377 \times 10^{5}$ on the right axis shows that curvature at the hotspots is already drastically altered for any $\beta > 0$.}
    \label{fig:L-curve}
\end{figure}

\begin{figure}[htbp]
    \centering
    \begin{minipage}{0.49\columnwidth}
        \centering
        \includegraphics[width=\linewidth]{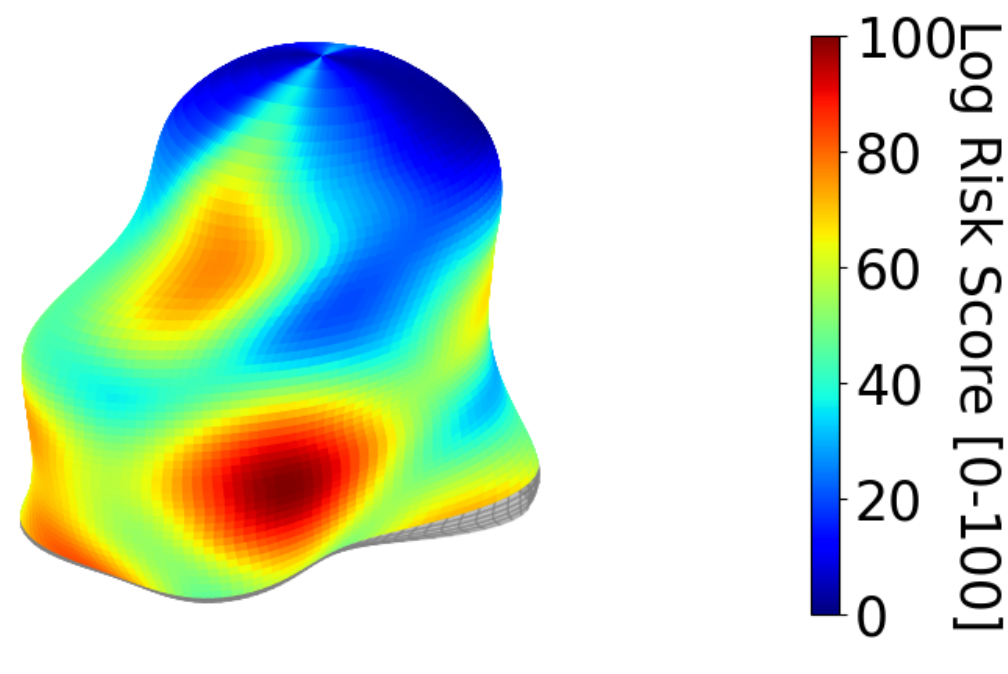}
        \vspace{0.25cm} (a)
    \end{minipage}\hfill
    \begin{minipage}{0.49\columnwidth}
        \centering
        \includegraphics[width=\linewidth]{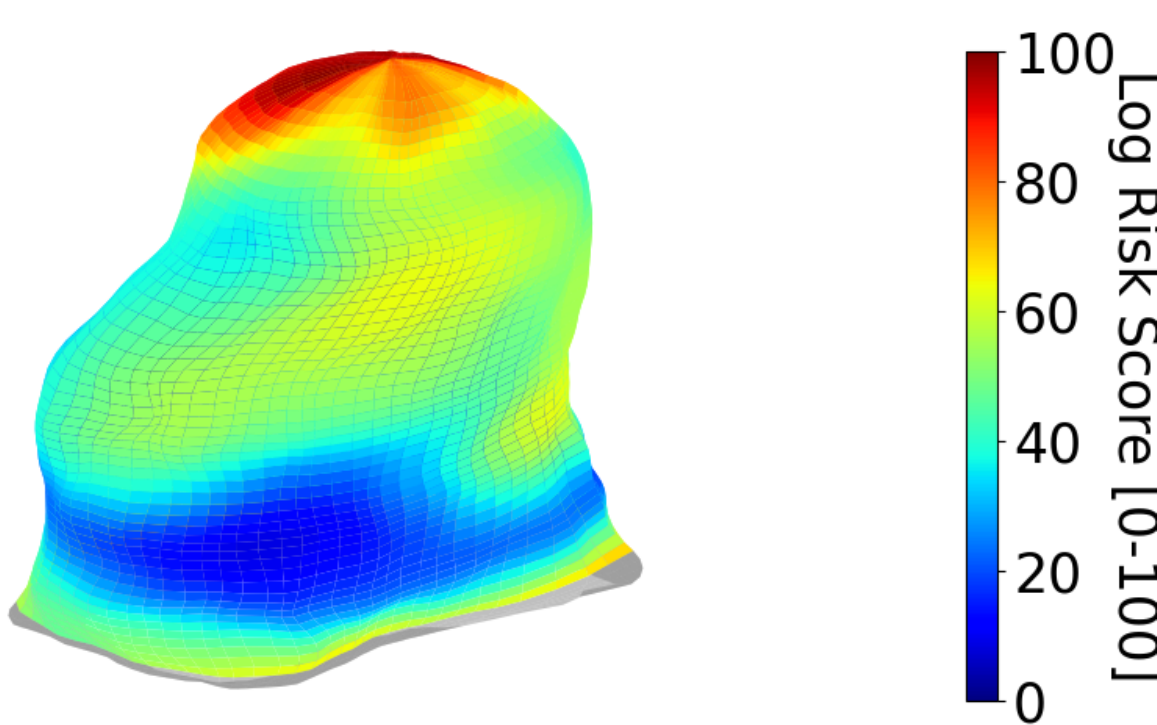}
        \vspace{0.5cm} (b)
    \end{minipage}
    
    \vspace{0.3cm}
    
    \begin{minipage}{0.49\columnwidth} 
        \centering
        \includegraphics[width=\linewidth]{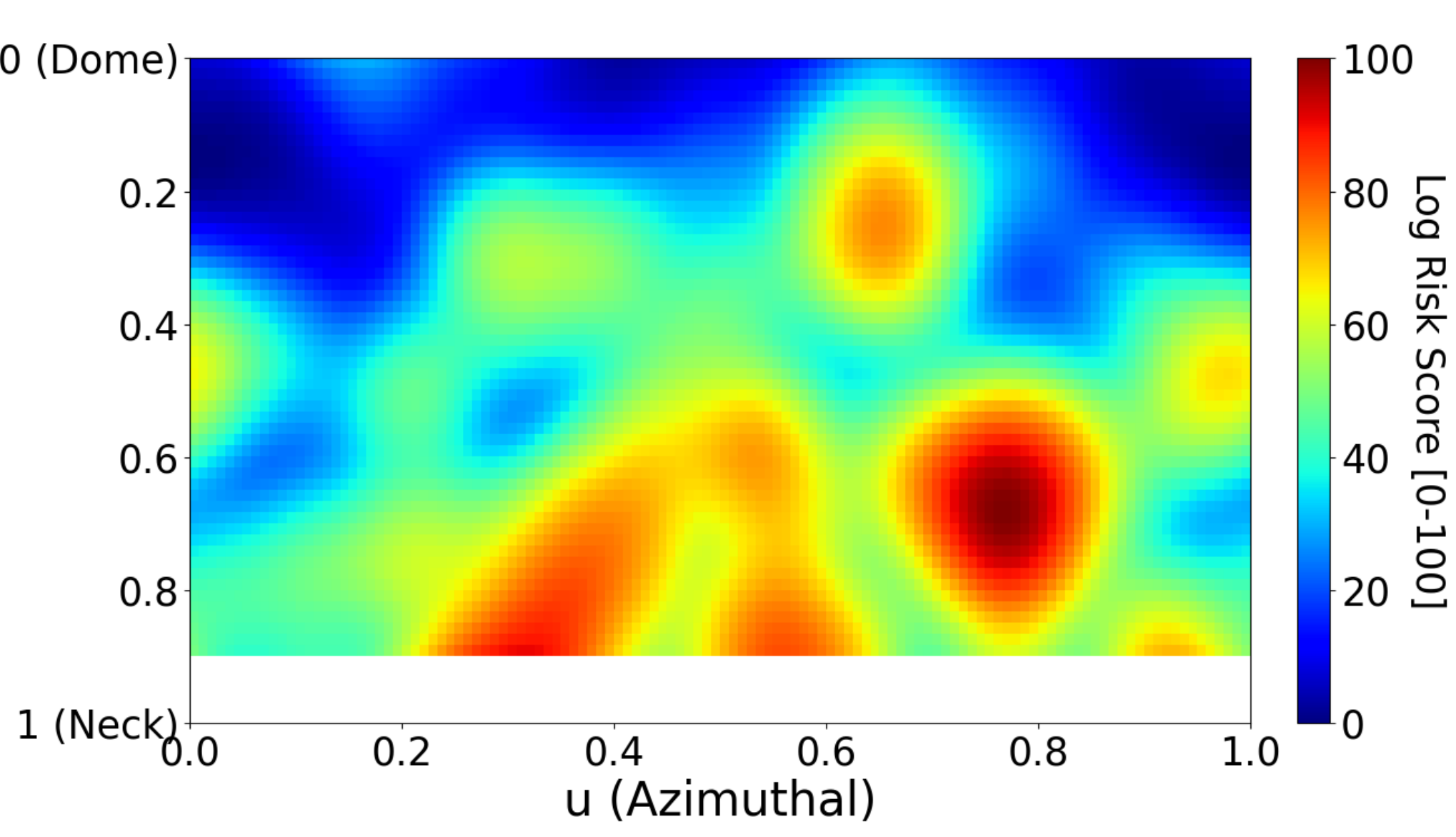}
        \vspace{0.2cm} (c)
    \end{minipage}\hfill
    \begin{minipage}{0.49\columnwidth}
        \centering
        \includegraphics[width=\linewidth]{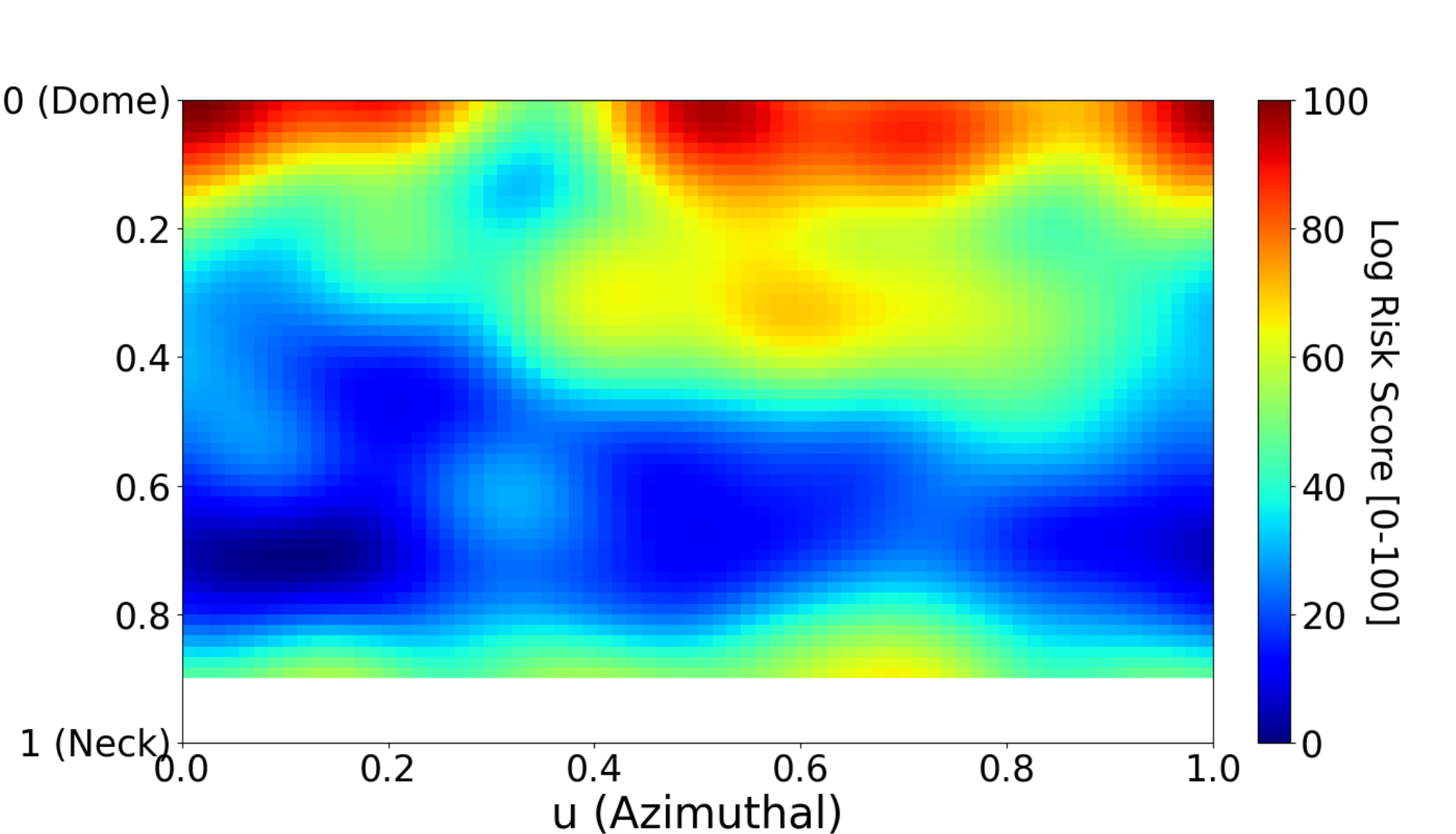}
        \vspace{0.2cm} (d)
    \end{minipage}
    \caption{Comparison of rupture risk maps between the conventional L-curve method and the proposed PINN framework (log risk score, 0--100). (a)~L-curve-based risk map on the 3D surface: hotspots (red/orange) are scattered across the mid-body and lower regions without a coherent pattern. (b)~PINN-based risk map on the 3D surface: the highest risk (saturated red) concentrates at the dome apex with a smooth gradient toward the neck (blue). (c)~L-curve-based risk map in the parametric domain ($u$--$v$): multiple isolated high-risk patches appear throughout the domain. (d)~PINN-based risk map in the parametric domain: a continuous high-risk band localizes near the dome ($v \approx 0$) and decays monotonically toward the neck ($v \approx 1$).} 
    \label{fig:risk_comparison}
\end{figure}

The constrained fitting framework employs the L-curve method~\citep{hansen1992analysis} to determine the optimal smoothing factor $\beta$. In log--log space, plotting the residual norm against the smoothness norm traces an L-shaped trajectory (Figure~\ref{fig:L-curve}a): the upper-left branch corresponds to over-smoothing (large $\beta$, high residual), while the lower-right branch corresponds to under-smoothing (small $\beta$, poor regularity). The elbow of this curve marks the optimal trade-off, yielding $\beta_{\mathrm{opt}} = 3.5 \times 10^{-7}$. However, because $\beta$ globally alters the surface curvature, any choice of smoothing factor systematically distorts the subsequent risk analysis. Figure~\ref{fig:L-curve}b quantifies these geometric and curvature deviations using the following metrics:
\begin{equation}
    \mathrm{SD}_{\mathrm{RMS}} = \sqrt{\frac{\sum (S(u,v)-Q)^2}{n_{Q}}}
\end{equation}
 
\begin{equation}
    C_{\mathrm{RMS}} = \sqrt{\frac{\sum_{i=1}^{10} (H(u_{i},v_{i},\beta)-H(u_{i},v_{i},0))^2}{10}}
\end{equation}
where $\mathrm{SD}_{\mathrm{RMS}}$ is the root-mean-square spatial deviation between the fitted surface and the measured points, $C_{\mathrm{RMS}}$ is the root-mean-square curvature error evaluated at the ten highest-risk hotspots, $n_{Q}$ is the number of measured points, and $(u_{i},v_{i})$ are the parametric coordinates of the ten highest-risk hotspots identified at $\beta=0$. Although the L-curve method identifies a mathematically optimal $\beta$, the conventional smoothing approach inevitably suppresses high-curvature regions. Specifically, $C_{\mathrm{RMS}}$ exhibits a persistent offset of $3.1377 \times 10^5$ with only marginal variation (order of $1$--$7$) across the entire parameter range (Figure~\ref{fig:L-curve}b, right axis). This indicates that the curvature at the hotspots is already drastically altered relative to the unsmoothed reference ($\beta = 0$) as soon as any smoothing is applied; no choice of $\beta$ can recover the original high-curvature features, leaving no viable window for parameter tuning. The resulting risk maps (Figure~\ref{fig:risk_comparison}) confirm this limitation.

As a result of L-curve optimization (Figure~\ref{fig:risk_comparison}a and c), the risk hotspots are scattered without a clear pattern. This contradicts a well-established finding from extensive clinical statistical literature and various computational analyses, which consistently identifies the highest rupture risk hotspot in ISAs at the apex of the dome~\citep{dhar2008morphology, cebral2010hemodynamics, baek2010clinical,torii2008fluid,humphrey2000structure,tawk2023pediatric}.
In contrast, the PINN-based risk map (Figure~\ref{fig:risk_comparison}b and d) successfully filters out these artifacts, revealing a clean and localized risk profile. The physics-informed regularization ensures that the reconstructed surface is consistent with biomechanical principles, producing a more clinically relevant assessment of rupture risk~\citep{brisman2006cerebral, inagawa2009risk, bacigaluppi2014factors, matsukawa2013morphological, cui2025prevalence}.

From a neurosurgical perspective, ruptured ISAs predominantly exhibit rupture points at or near the dome apex, where the arterial wall is maximally attenuated~\citep{dhar2008morphology, brisman2006cerebral}. The PINN-based risk map (Figure~\ref{fig:risk_comparison}b) concentrates the highest risk scores at the dome apex, consistent with this well-established finding. By contrast, the spatially scattered hotspot pattern from the L-curve method (Figure~\ref{fig:risk_comparison}a) does not correspond to any recognizable pattern of wall vulnerability observed during microsurgical clipping or endovascular coiling. If used for preoperative planning, such dispersed distributions could mislead the surgical team by falsely indicating multiple regions of equal rupture susceptibility---a scenario rarely encountered in clinical practice.

\section{Discussion}

This study addressed the problem that CTA/MRA-derived surface reconstructions of intracranial saccular aneurysms (ISAs) inherit hardware-level imaging noise, which corrupts curvature-based rupture risk assessment. The proposed PINN framework selectively suppresses imaging artifacts while preserving genuine geometric features by enforcing Laplace membrane equilibrium as a physics-informed constraint.
The results demonstrate that the physics-informed framework produces more reliable risk maps than conventional constrained fitting with smoothing. The following subsections discuss the methodological contributions, inherent limitations, and future extensions of the proposed framework.

\subsection{Contribution}
As discussed in Section 1, CTA/MRA images possess hardware-level constraints, such as a poor signal-to-noise ratio (SNR) and partial volume effects, which induce imaging artifacts. Notably, even conventional catheter-based digital subtraction angiography (DSA), often regarded as the reference standard, is not immune to noise and inherently lacks the capacity to distinguish clinically significant lesions---such as blebs or rupture points---from imaging artifacts. While conventional smoothing filters at each representation level—such as anisotropic diffusion at the pixel intensity level \citep{gerig1992nonlinear}, curvature-based filters at the discrete mesh vertex level \citep{taubin1995signal, cebral2010hemodynamics}, and Tikhonov regularizations at the parametric control point level \citep{hansen1992analysis, styner2006framework}—effectively mitigate these constraints for statistical-data-based clinical judgments, they introduce critical challenges in computational simulations. In rupture risk assessment, it is paramount to differentiate artifact-compensating smoothing from true, fine-geometric features; for instance, high-risk hotspots like blebs—which are extensively documented as critical regions of localized pressure concentration and hemodynamic deterioration \citep{ashkezari2021blebs, cebral2010bleb, suzuki2026hemodynamic}— are often numerically indistinguishable from non-physical artifacts.
To overcome this fundamental limitation, the present study proposes a physics-informed neural network (PINN) framework. By incorporating the governing equations of ISA mechanical behavior into the physics loss term, the iterative optimization constraint robustly drives the system toward a realistic Laplace membrane equilibrium, thereby successfully filtering out imaging artifacts without compromising genuine geometric anomalies. This physics-informed approach ensures that the reconstructed surface is strictly consistent with biomechanical principles, resulting in a more clinically relevant and reliable assessment of rupture risk~\citep{hughes1998variational, bendsoe2003topology, boo2026isogeometric}. The key contributions of this framework are (i) embedding the variational equilibrium condition ($\delta \Pi = 0$) directly into the geometric reconstruction pipeline as a physics-informed loss, (ii) developing a Manifold-Consistent CNN ansatz that preserves the topological structure of closed vascular surfaces, and (iii) establishing a physics-driven criterion---Laplace membrane equilibrium rather than the mathematical L-curve---to discriminate genuine pathological features from imaging artifacts.

From a clinical perspective, the identification of biomechanical weak points (Figure~\ref{fig:risk_comparison}) offers several practical implications for neurosurgical practice. First, during follow-up surveillance of unruptured aneurysms, clinicians currently rely primarily on CTA/MRA to monitor size changes; incorporating weak point pattern evolution could provide an additional dimension for longitudinal risk assessment. Second, changes in weak point distribution over time may serve as a predictive indicator to differentiate patients requiring preventive intervention from those who can be safely monitored. Third, during preventive procedures such as clip ligation or coil embolization, prior knowledge of weak point locations enables the surgical team to exercise heightened caution near vulnerable regions, potentially reducing procedural complication rates.

This work bridges geometric modeling and mechanical analysis in computational biomedicine. Rather than treating them as separate sequential stages, the stationarity of the total potential energy functional ($\delta \Pi = 0$) is used to define the physical state of the geometry itself: the final reconstructed shape is the equilibrium configuration at the end of deformation, providing a physics-driven denoising mechanism. From a clinical standpoint, this integration enables patient-specific rupture risk mapping directly from routine CTA/MRA data, without requiring invasive measurements or additional imaging modalities. By converting standard neuroimaging into a biomechanically grounded risk profile, the framework has the potential to support clinical decision-making in the surveillance and treatment planning of unruptured ISAs.

\subsection{Limitations}

This research is not without limitations.
First, the aneurysm wall is modeled as an isotropic elastostatic membrane. Although mature aneurysms tend toward isotropic stiffness profiles due to elastin degradation, this assumption neglects wall thickness variation, anisotropic collagen fiber architecture, and hemodynamic loading~\citep{shah1999finite, cebral2005characterization, busch2011statistical}. Incorporating fluid--structure interaction (FSI) or computational fluid dynamics (CFD) under pulsatile blood flow would provide additional indicators such as wall shear stress gradients for patient-specific diagnosis~\citep{berg2017does, valen2018realworld}.

Second, the current validation relies on qualitative agreement with neurosurgical experience rather than direct comparison against confirmed rupture site data. The side wall aneurysm examined in this study presents an inherent validation challenge: the majority of ruptured intracranial aneurysms occur at bifurcation sites, and side wall aneurysm ruptures are relatively uncommon. Moreover, when side wall aneurysms do rupture, endovascular coil embolization is generally preferred over surgical clipping, rendering direct intraoperative visual confirmation of the rupture point infeasible in most cases.

\subsection{Future extensions}
The current framework faces parameterization complexities when applied to multilobed aneurysms, where overlapping folds exceed the mapping capabilities of single-patch manifold-consistent padding~\citep{thani2026physics, zhu2025graph}. Future work will focus on multi-patch NURBS integration and dynamic loading conditions to improve diagnostic fidelity. Extending the framework to bifurcation aneurysms---where surgical clipping is more frequently performed---could enable direct comparison between predicted weak points and intraoperatively confirmed rupture sites, thereby providing quantitative clinical validation. In addition, longitudinal tracking of weak point pattern evolution across serial imaging may offer a predictive biomarker for aneurysm growth and rupture risk stratification.

\section{Conclusion}

A physics-informed neural network (PINN) framework has been developed for Laplace equilibrium-driven surface reconstruction aimed at patient-specific rupture risk assessment of intracranial saccular aneurysms (ISAs). Building upon the classical Laplace membrane equilibrium ($\kappa_1 T_1 + \kappa_2 T_2 = P$)~\citep{humphrey1996use}, the proposed framework integrates this governing relation into an image-based geometric reconstruction pipeline via deep learning---a connection that has remained unexplored. The key contributions are: (i)~embedding the variational equilibrium condition ($\delta \Pi = 0$) directly into the reconstruction pipeline as a physics-informed loss function, (ii)~developing a Manifold-Consistent CNN ansatz that preserves the topological structure of closed vascular surfaces, and (iii)~establishing a physics-driven criterion---Laplace membrane equilibrium rather than the mathematical L-curve---to discriminate genuine pathological features from imaging artifacts. By enforcing this physics-driven constraint, the framework selectively suppresses hardware-induced imaging artifacts in CTA/MRA data without compromising genuine high-curvature features such as blebs, yielding a clinically consistent rupture risk profile with risk concentration at the dome apex.
While the current single-patch isotropic membrane formulation is limited when handling multilobed morphologies (e.g., cactus-shaped aneurysms), it provides a foundation for patient-specific modeling. Incorporating multi-patch NURBS, fluid--structure interaction (FSI), and pulsatile flow conditions represents a clear path for future work. From a clinical standpoint, the framework enables patient-specific rupture risk mapping directly from routine CTA/MRA data without requiring invasive measurements, offering a computational biomedicine tool to support clinical decision-making in the surveillance and treatment planning of unruptured ISAs.

\section*{Data Availability}

The data that support the findings of this study are available from the corresponding author upon reasonable request. The data are not publicly available due to privacy or ethical restrictions.

\section*{Acknowledgments}
\begin{itemize}
    \item This work was supported by the National Research Foundation of Korea(NRF) grant funded by the Korea government(MSIT) (No.RS-2026-25501401).
    \item This research was funded by the `Changwon National University - Samsung Changwon Hospital joint Collaboration Research Support Project' in 2025.
    \item This research was funded by the `Samsung Changwon Hospital - Changwon National University Joint Collaboration Research Support Project' in 2025.
    % \item This study was conducted as part of the Glocal University Project, supported by the RISE (Regional Innovation System \& Education) program funded by the Ministry of Education.
    % \item This research was supported by the Regional Innovation System \& Education(RISE) program through the RISE Center, Gyeongsangnam-do, funded by the Ministry of Education(MOE) and the Gyeongsangnam-do Provincial Government, Republic of Korea.(2026-RISE-16-002).
\end{itemize}

\section*{Declaration of Generative AI and AI-assisted technologies in the writing process}
During the preparation of this work, the authors used OpenAI's ChatGPT in order to improve clarity and readability. After using this tool, the authors reviewed and edited the content as needed and take full responsibility for the content of the published article.

\section*{Declaration of competing interest}
The authors declare that they have no known competing financial interests or personal relationships that could have appeared to influence the work reported in this paper.

\section*{CRediT Author Contributions}
\textbf{Hyomin Ryu}: Writing -- original draft, Software, Validation, Formal analysis, Investigation, Data curation, Visualization.
\textbf{Seung Hwan Kim}: Resources, Writing -- review \& editing, Funding acquisition, Ethics approval (IRB).
\textbf{Jaemin Kim}: Conceptualization, Methodology,  Data curation, Writing -- review \& editing, Visualization, Supervision, Project administration, Funding acquisition.
\section*{Appendix A. Derivation of physical grid dimensions}
\renewcommand{\theequation}{A\arabic{equation}}
\setcounter{equation}{0}
\renewcommand{\thefigure}{A\arabic{figure}}
\setcounter{figure}{0}
To establish a mathematically rigorous foundation for the geometric reconstruction, the spatial degrees of freedom for the baseline control point grid were determined by evaluating the directional resolutions of the imaging hardware within the curvilinear domain of the aneurysm. The source dataset was acquired using a 3.0-Tesla Magnetic Resonance Imaging (MRI) system via a 3D Time-of-Flight (TOF) Brain MRA sequence, yielding an in-plane spatial resolution ($\Delta x_{img} = 0.390625\text{ mm}$) and a slice thickness ($\Delta z_{img} = 0.600000\text{ mm}$).

The target intracranial saccular aneurysm, featuring a clinical height ($L$) of $6.29\text{ mm}$ and a neck diameter ($D$) of $4.53\text{ mm}$ (Figure~\ref{fig:slicer_workflow}C), was modeled as an ellipsoidal hemisphere~\citep{ma2004three, dhar2008morphology}. The macroscopic physical boundaries were quantified into the azimuthal circumference ($c_\theta$) at the base and the polar arc length ($s_\phi$) from the apex to the base. While $c_\theta$ is directly computed from the neck diameter ($c_\theta = \pi D \approx 14.23\text{ mm}$), $s_\phi$ represents a quarter-arc of an ellipse with semi-axes $a = D/2 = 2.265\text{ mm}$ and $b = L = 6.29\text{ mm}$. Utilizing Ramanujan's first asymptotic approximation for ellipse perimeters, these macro-scale geometric lengths were formalized as follows~\citep{zhang2007patient, fessler2010model}:
\begin{equation}
c_\theta = \pi D \approx 14.23\text{ mm}
\end{equation}
\begin{equation}
s_\phi \approx \frac{\pi}{4} \left[ 3(a+b) - \sqrt{(3a+b)(a+3b)} \right] \approx 7.10\text{ mm}
\end{equation}

Concurrently, because the imaging hardware possesses anisotropic resolution thresholds ($\Delta x_{img} \neq \Delta z_{img}$), a generalized maximum arc-direction resolution unit ($\Delta s_{res}$) was established to map the Nyquist-Shannon sampling criterion \citep{unser2000sampling} consistently onto the curved manifold. By treating the hardware thresholds as semi-axes ($a_{res} = \Delta x_{img}$, $b_{res} = \Delta z_{img}$), the representative quarter-arc resolution unit was derived via the same Ramanujan formulation:

\begin{equation}
\Delta s_{res} \approx \frac{\pi}{4} \left[ 3(a_{res}+b_{res}) - \sqrt{(3a_{res}+b_{res})(a_{res}+3b_{res})} \right] \approx 0.79\text{ mm}
\end{equation}

Mapping this physical hardware limit directly to the structural dimensions of the aneurysm yields the theoretical maximum physical grid dimensions ($n_{\theta, phys}, n_{\phi, phys}$) required to fully capture the geometry without aliasing or information loss:
\begin{equation}
n_{\theta, phys} = \left\lfloor \frac{c_\theta}{\Delta s_{res}} \right\rfloor = \left\lfloor \frac{14.23}{0.79} \right\rfloor = 18
\end{equation}
\begin{equation}
n_{\phi, phys} = \left\lfloor \frac{s_\phi}{\Delta s_{res}} \right\rfloor = \left\lfloor \frac{7.10}{0.79} \right\rfloor = 9
\end{equation}

The strictly derived physical grid of $18 \times 9$ was minorly adjusted to a final optimized control point network of $16 \times 8$ ($n_\theta = 16$, $n_\phi = 8$) to comply with the architectural constraints of Manifold-consistent CNN. Synchronizing the initial degrees of freedom to a $16 \times 8$ ensures strict mathematical consistency during the spatial padding operations and subsequent dyadic upsampling steps ($16 \times 8 \rightarrow 64 \times 32$) within the convolutional trial function pipeline~\citep{ronneberger2015unet, long2015fully}. Functioning as an embedded low-pass filter, this strategically bounded grid filters out voxel-induced staircase artifacts ($\le 0.6\text{ mm}$) and coordinate-transformation noise (Figure~\ref{fig:extract_ptc}), while safely retaining the essential macro-geometric features and structural traces of pathological blebs for the conventional constrained fitting detailed in Section 2.3.

\begin{figure}[H]
  \centering
  \includegraphics[width=\linewidth]{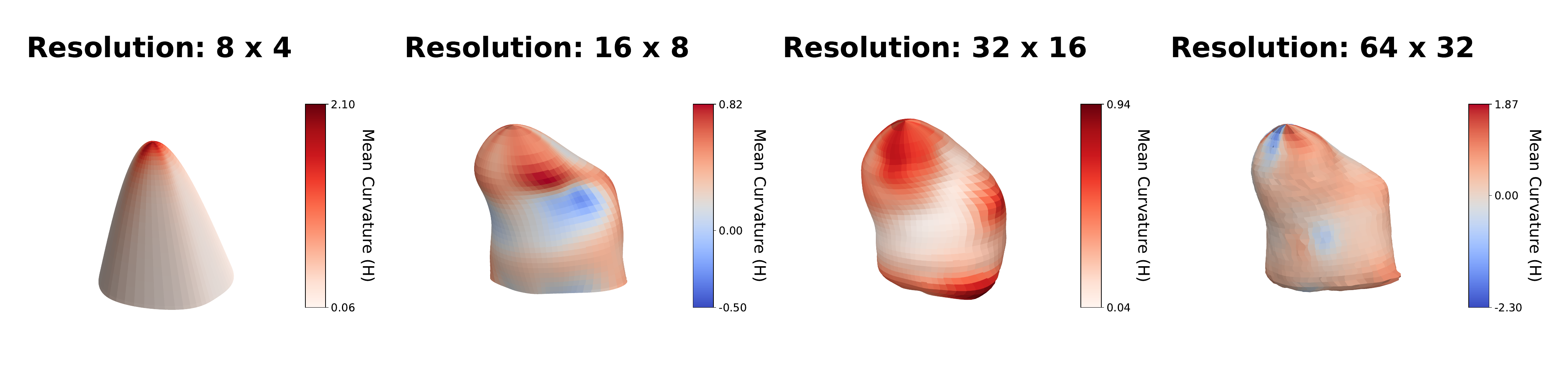}
  \caption{The surfaces for each resolution. Intuitively, resolutions below 16$\times$8 are insufficient to represent the geometry, while higher resolutions begin to capture the aneurysm's shape.}
  \label{fig:resolution_comparison}
\end{figure}

\section*{Appendix B. Mathematical formulation}
\renewcommand{\theequation}{B\arabic{equation}}
\setcounter{equation}{0}

\subsection*{B.1 System potential and work done}
The analytical model is based on the following fundamental assumptions:
\begin{enumerate}
    \item $P = \kappa_1 T_1 + \kappa_2 T_2$ \quad (Laplace equation)~\citep{kim2023mechanosensitive, kim2023modelingmechanosensitivecollectivemigration}
    \item Steady-state condition
    \item Membrane mechanics
\end{enumerate}

Let $d$ denote the displacement field. The system potential ($\Pi$) is defined as the difference between the internal strain energy ($U$) and the external work done ($W$):
\begin{equation}
    \Pi = U - W
\end{equation}
where the strain energy $U$ consists of bending and stretching components, i.e., $U = U_{\text{bending}} + U_{\text{stretch}}$.

Based on the Humphrey Growth and Remodeling (G\&R) theory, the external work done is given by:
\begin{equation}
    \int_V P \, dV = \int_V (\kappa_1 T_1 + \kappa_2 T_2) \, dV = T_h t \int_S (\kappa_1 + \kappa_2) \, dA
\end{equation}
assuming uniform tension $T_1 \approx T_2 = T_h$ (constant) and the volume element $dV \approx t \, dA$ (where $t$ is thickness and $dA$ is the surface area element).

\subsection*{B.2 Incremental linearization}
Under the steady-state assumption and transition from nonlinear hyperelastic to linear elastic behaviors (incremental linearization), the membrane bending energy is assumed to be negligible:
\begin{equation}
    U_{\text{bending}} \approx \frac{E t^3}{1 - \nu^2} \approx 0 \quad \text{(Membrane)}
\end{equation}

The stretching energy is fundamentally expressed as:
\begin{equation}
    U_{\text{stretch}} = \frac{1}{2} E \epsilon^2 V
\end{equation}
where the strain components with respect to the coordinates are approximated as $\epsilon_u \approx \left\| \frac{\partial d}{\partial u} \right\|$ and $\epsilon_v \approx \left\| \frac{\partial d}{\partial v} \right\|$. 

The directional strain $\epsilon_\theta$ can be expanded as:
\begin{equation}
    \epsilon_\theta \approx \frac{\partial d}{\partial \theta} = \cos\theta \frac{\partial d}{\partial u} + \sin\theta \frac{\partial d}{\partial v}
\end{equation}
Integrating the square of the directional strain over a full cycle ($2\pi$) yields:
\begin{equation}
    \int_0^{2\pi} \epsilon_\theta^2 \, d\theta = \pi \left( \left\| \frac{\partial d}{\partial u} \right\|^2 + \left\| \frac{\partial d}{\partial v} \right\|^2 \right)
\end{equation}

By defining the stiffness parameter $C = \pi t E$, the primary stretching energy can be effectively approximated as:
\begin{equation}
    U_{\text{stretch}} \approx \frac{1}{2} C \iint_S \left( \left\| \frac{\partial d}{\partial u} \right\|^2 + \left\| \frac{\partial d}{\partial v} \right\|^2 \right) du \, dv
\end{equation}

\subsection*{B.3 Variational principle and compliance}

Applying the variational principle to minimize the system potential ($\delta \Pi = 0$), the variation with respect to the displacement field $d$ is obtained as:
\begin{equation}
    \therefore \frac{\partial \Pi}{\partial d} = \lambda \iint_S \left( \left\| \frac{\partial d}{\partial u} \right\| + \left\| \frac{\partial d}{\partial v} \right\| \right) du \, dv - \frac{\partial \left( \int_S (\kappa_1 + \kappa_2) \, dA \right)}{\partial d}
\end{equation}
where the parameter $\lambda$ is given by:
\begin{equation}
    \lambda = \frac{C}{t T_h} \approx 500
\end{equation}

\subsection*{B.4 Tangential suppression}
A tangential displacement suppression penalty is introduced. For the baseline B-spline surface $\mathbf{S}_{base}(u, v)$ from $\mathbf{P}_{base}$, the covariant basis vectors $\mathbf{g}_u$ and $\mathbf{g}_v$ are computed via finite differences:
\begin{equation}
    \mathbf{g}_u = \frac{\partial \mathbf{S}_{base}}{\partial u}, \quad \mathbf{g}_v = \frac{\partial \mathbf{S}_{base}}{\partial v}
\end{equation}
The local unit normal vector $\hat{\mathbf{n}}$ is subsequently defined via the cross product: 
\begin{equation}\label{eq:normal_vector}
    \hat{\mathbf{n}} = \frac{\mathbf{g}_u \times \mathbf{g}_v}{\|\mathbf{g}_u \times \mathbf{g}_v\| + \epsilon_{norm}}
\end{equation}
where $\epsilon_{norm} = 10^{-12}$ is a sufficiently small stabilization constant. Under this local coordinate system, the predicted control point displacement field $\mathbf{d}$ is orthogonally decomposed into normal and tangential components. The tangential component $\mathbf{d}_{tangent}$ is extracted by subtracting the projected normal displacement from the total displacement:
\begin{equation}\label{eq:tangent_decomp}
    \mathbf{d}_{tangent} = \mathbf{d} - \left( \mathbf{d} \cdot \hat{\mathbf{n}} \right) \hat{\mathbf{n}}
\end{equation} 

\section*{Appendix C. Convergence analysis}
\renewcommand{\theequation}{C\arabic{equation}}
\setcounter{equation}{0}
\renewcommand{\thefigure}{C\arabic{figure}}
\setcounter{figure}{0}

\begin{figure}[H]
    \centering
    \includegraphics[width=1\linewidth, height=0.55\linewidth]{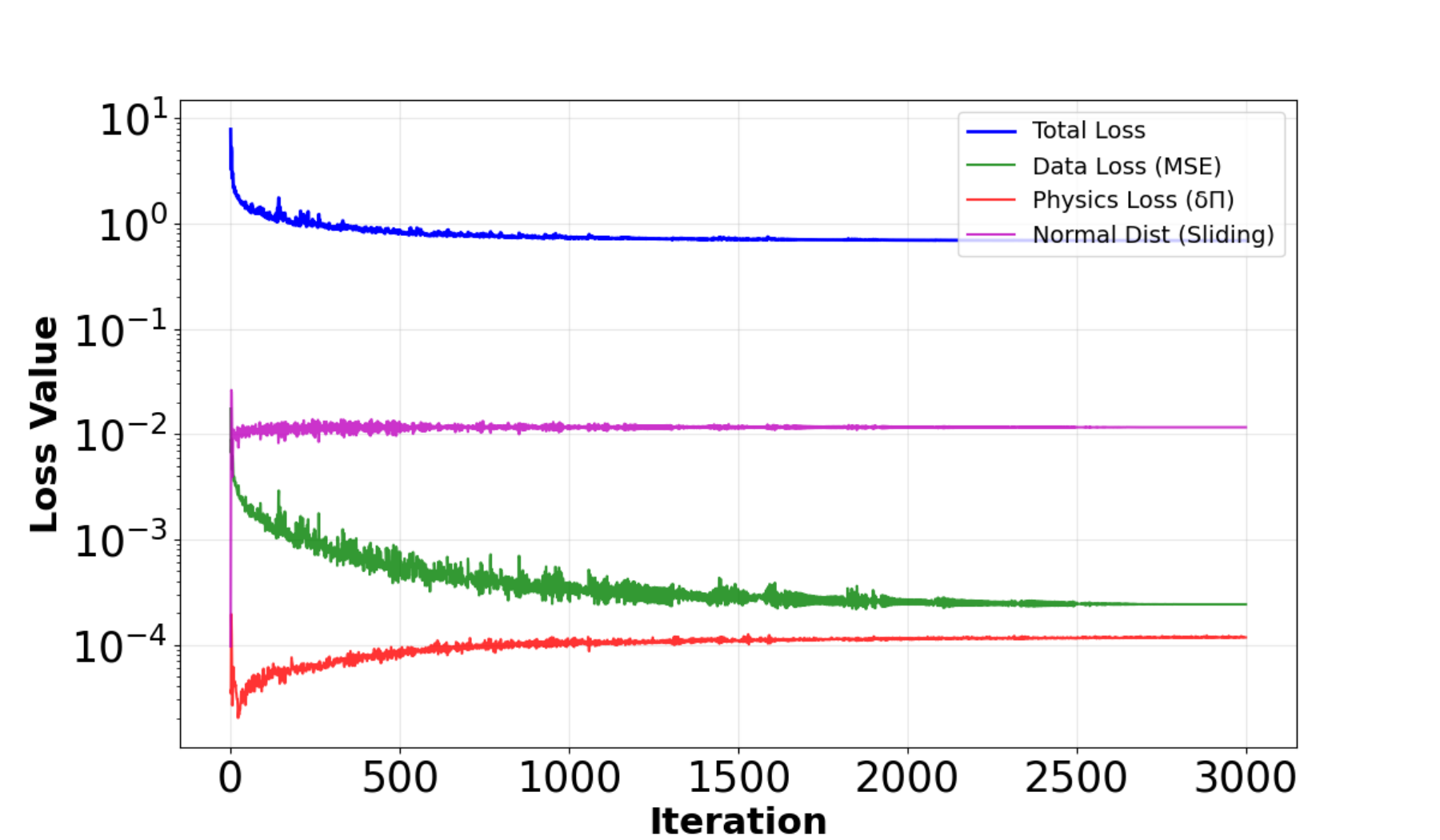}
    \caption{Loss curves. The Data Loss and Physics Loss exhibit intense waving and competition, with the amplitude of each curve gradually mitigating as they converge.}
    \label{fig:loss_and_recon}
\end{figure}

This framework was implemented using the CNN-Ansatz architecture initialized with the macro geometry, and optimized via advanced gradient descent strategies~\citep{loshchilov2016sgdr, loshchilov2017decoupled}. As illustrated in Figure~\ref{fig:loss_and_recon}, let us first examine the convergence of the loss curves. Both the data-driven MSE and the variational equilibrium loss successfully converged around 2,000 iterations. Notably, the data fidelity loss exhibited severe oscillations prior to 2,000 iterations. This drastic fluctuation demonstrates the intense competitive interplay between the data loss, which attempts to overfit imaging artifacts, and the physics loss, which strictly enforces mechanical regularization to reject non-physical geometries. Crucially, the data loss term serves a dual purpose: it not only maintains the overall macro-shape against unconstrained deformation, but also resolves the inherent non-uniqueness of the variational equilibrium by anchoring the surface to the measured points. Following this dynamic equilibration phase, the reconstructed surface successfully approximates a state of static equilibrium between the internal membrane tension and external blood pressure.

\clearpage

\bibliographystyle{elsarticle-harv}
\bibliography{references}

\end{document}